\documentclass[11pt,twoside]{article}
\usepackage[T1]{fontenc}
\usepackage{lmodern}
\usepackage[a4paper,margin=24mm,headheight=28pt,headsep=10pt]{geometry}
\usepackage{amsmath,amssymb,amsthm,mathtools}
\usepackage{microtype}
\usepackage{xcolor}
\newif\ifrevision
\revisiontrue

\usepackage[numbers,sort&compress]{natbib}
\usepackage[colorlinks=true,linkcolor=blue,citecolor=blue,urlcolor=blue]{hyperref}
\theoremstyle{plain}
\newtheorem{theorem}{Theorem}[section]
\newtheoremstyle{uprightremark}%
  {\topsep}{\topsep}{\normalfont}{}%
  {\bfseries}{.}{0.5em}%
  {\thmname{#1}\thmnumber{ #2}\thmnote{ \textbf{(#3)}}}
\theoremstyle{uprightremark}
\newtheorem{remark}{Remark}[section]
\newcommand{\R}{\mathbb{R}}

\newcommand{\Om}{\Omega}
\usepackage{tikz}
\usetikzlibrary{arrows.meta}

\usepackage{caption}
\definecolor{KSblue}{RGB}{40,105,170}
\definecolor{KSorange}{RGB}{220,110,35}
\definecolor{KSteal}{RGB}{30,145,140}

\usepackage{hyperref}
\hypersetup{
	colorlinks=true,
	citecolor=red,
	linkcolor=blue
}
\let\originalfootnote\footnote
\renewcommand{\footnote}[1]{%
	\begingroup
	\hypersetup{linkcolor=KSorange}%
	\originalfootnote{#1}%
	\endgroup
}

\usepackage{mathrsfs,bm,needspace}
\newcommand{\into}{\int_\Omega}
\newcommand{\Tmax}{T_{\max}}

\numberwithin{equation}{section}
\allowdisplaybreaks[1]
\theoremstyle{plain}
\newtheorem{lem}{Lemma}[section]
\newtheorem{prop}{Proposition}[section]
\theoremstyle{uprightremark}
\newtheorem{rem}{Remark}[section]
\newcommand{\RunningTitle}{Optimal mass threshold for the 2D Keller-Segel-Navier-Stokes system}
\usepackage{fancyhdr}
\newcommand{\PaperHeaderSetup}{%
  \fancyhf{}%
  \fancyhead[RO]{\footnotesize\nouppercase{\RunningTitle}}%
  \fancyhead[LE]{\footnotesize\nouppercase{L. Du, Y. Wang, and Y. Zhang}}%
  \fancyfoot[C]{\thepage}%
  \renewcommand{\headrulewidth}{0.4pt}%
}
\fancypagestyle{paperstyle}{\PaperHeaderSetup}
\fancypagestyle{plain}{%
  \fancyhf{}%
  \fancyfoot[C]{\thepage}%
  \renewcommand{\headrulewidth}{0pt}%
}
\title{On the optimal mass threshold for the two-dimensional Keller-Segel system coupled to the Navier-Stokes equations}
\author{
	Lili Du$^1$\footnote{dulili@scu.edu.cn}, \,\,Yulan Wang$^2$\footnote{wangyulan-math@163.com}, \,\,Yichen Zhang$^1$\footnote{yczhang@stu.scu.edu.cn}\\[2mm]
	\small$^1$Department of Mathematics, Sichuan University, Chengdu 610065, China\\
	\small$^2$School of Science, Xihua University, Chengdu 610039, China
	}
\date{}

\begin{document}
\maketitle
\thispagestyle{plain}
\begin{abstract}
We investigate the mass threshold for the two-dimensional Keller-Segel system coupled with the incompressible Navier-Stokes equations. In the absence of fluid interaction, the classical Keller-Segel system admits a well-known $4\pi$-critical-mass dichotomy: for every suitably regular initial datum of total mass \(m<4\pi\), the associated no-flux initial-boundary value problem admits a global bounded classical solution, whereas for \(m>4\pi\), there exist smooth initial configurations of mass  $m$ such that the associated solution blows up in finite or infinite time. {For the fluid-coupled system, Winkler ({SIAM J. Math. Anal.}, 52:2041-2080, 2020) first established global generalized solvability under the condition  \(m<2\pi\) and posed an open question whether the fluid action may lower the critical blow-up mass to some value lying between  $2\pi$ and $4\pi$. Based on numerical evidence, Bonilla and Guti\'errez-Santacreu ({J. Numer. Math.}, 32:175-212, 2024) later conjectured that the critical value is still \(4\pi\).

In the present work, we resolve this open problem and provide a rigorous mathematical proof to the conjecture, which indicates that} the presence of the fluid does not alter the critical mass threshold, and the entire classical subcritical regime remains intact---every admissible initial datum with total mass \(m<4\pi\) generates a unique global classical solution that is uniformly bounded in time.

Our proof relies on the construction of a new quasi-energy inequality. The key underlying observation is that combining relative isoperimetric estimates with rearrangement techniques yields a Gagliardo-Nirenberg inequality on bounded planar domains with an explicit leading constant.\\[2mm]
\noindent\textbf{2020 Mathematics Subject Classification.} 92C17, 35Q35, 35Q92, 35A09.\par
\noindent\textbf{Keywords.} Keller-Segel-Navier-Stokes system, critical mass, global existence, boundedness.
\end{abstract}


\section{Introduction}

The seminal Keller-Segel model~\cite{KellerSegel1970,KellerSegel1971} provides a fundamental description of biological aggregation through the competition between diffusion and attraction toward a self-produced chemical signal. Its broader biological significance and extensive mathematical theory are explored in the influential contributions and surveys~\cite{LevinSegel1985,HillenPainter2009,KiselevNazarovRyzhikYao2023,Horstmann2003,BellomoBellouquidTaoWinkler2015}. The question of a critical mass itself has deep biological roots: in a pioneering study, Nanjundiah~\cite[p.~102]{Nanjundiah1973} identified a density threshold for aggregation, noting that ``above a certain critical density, the cells can aggregate''. More specifically, Childress and Percus conjectured that ``if an aggregate contains more than a critical number of cells, the aggregate can collapse to infinite density (in finite time). For aggregates less than this critical size, such a collapse would then never occur'' (\cite{Childress1984,ChildressPercus1981}). These prescient biological insights and the subsequent landmark mathematical advances~\cite{NagaiSenbaYoshida1997,BlanchetDolbeaultPerthame2006,BlanchetCarrilloMasmoudi2008,HorstmannWang2001,Mizoguchi2020Cauchy,Nagai2001} underscore the enduring prominence of critical mass in Keller-Segel systems and mathematical biology.

For the classical two-dimensional Keller-Segel system in bounded domains, the celebrated work of Nagai, Senba, and Yoshida~\cite{NagaiSenbaYoshida1997} established global existence and uniform boundedness when the initial population mass \(m:=\int_\Omega n_0\) satisfies \(m<4\pi\), whereas  Horstmann and Wang \cite{HorstmannWang2001} demonstrated the sharpness of this mass bound through the existence of unbounded solutions for suitable supercritical data with \(m>4\pi\). Consequently, $4\pi$ is the critical mass for the two-dimensional Keller-Segel system in bounded domains.

Biological phenomena are sometimes modeled using Keller-Segel-fluid coupled systems because biological reality often requires accounting for the effects of the fluid on populations and chemical signals \cite{BellomoOutadaSolerTaoWinkler2022,BellomoBellouquidTaoWinkler2015}. Naturally, this raises a fundamental question: does fluid interaction affect the classical critical-mass threshold of the coupled system? In a pioneering contribution, Winkler~\cite{Winkler2020SIAM} first established global generalized solvability {for a two-dimensional Keller-Segel-Navier-Stokes system} under the explicit condition $m<2\pi$, and posed the following open problem \cite[p.~2045]{Winkler2020SIAM}:
\begin{quote}
	\itshape
	‘‘Interesting natural problems left open here are to clarify how far the explicit condition (1.7) indeed is optimal for the conclusion made in Theorem 1.1, $\cdots\cdots$. However, Theorem 1.1 leaves open the possibility that, $\cdots\cdots$, the fluid interaction mechanism in (1.3) might at least slightly affect the precise value of the critical blow-up mass from (1.2), in the sense of a reduction to a level between \(2\pi\) and \(4\pi\), if the initial data for the fluid velocity field satisfy appropriate constraints.’’\footnote{Here, in \cite{Winkler2020SIAM}, the notation (1.7) denotes the condition $\into n_0<2\pi$, while (1.3) and (1.2) refer to the Keller-Segel-Navier-Stokes system and its fluid-free counterpart, respectively.}
\end{quote}
{Remarkably, a striking numerical perspective on this question presented in \cite{BonillaGutierrezSantacreu2024} subsequently suggested that \(m<2\pi\) might not be optimal, as Bonilla and Guti\'errez-Santacreu pointed out in \cite[Abstract]{BonillaGutierrezSantacreu2024}:
\begin{quote}
	\itshape
	``Our findings show that the existence threshold value $2\pi$ encountered for the organism density mass may not be optimal and hence it is conjectured that the critical threshold value $4\pi$ may be inherited from the fluid-free Keller-Segel equations.''
\end{quote}
The present work decisively resolves this open problem posed by Winkler and gave an positive answer to this conjecture}, namely, every admissible initial datum with \(m<4\pi\) generates a unique global classical solution that remains uniformly bounded in time. Thus, incompressible fluid interaction does not lower the universal \(4\pi\) threshold for global boundedness.

\subsection{Critical mass for the classical Keller-Segel system}\label{sec:parabolic-ks}\label{sec:elliptic-comparison}

The critical mass for global existence and blow-up in the Keller-Segel system has received sustained attention, with numerous results accumulated over the past few decades.

For a bounded domain \(\Omega\subset\mathbb{R}^2\) with smooth boundary, the fully parabolic Keller-Segel system with homogeneous Neumann boundary conditions reads
\begin{equation}\label{eq:PP-Neumann}
	\begin{cases}
		n_t=\Delta n-\nabla\!\cdot(n\nabla c),&x\in\Omega,\ t>0,\\
		c_t=\Delta c-c+n,&x\in\Omega,\ t>0,\\
		\partial_\nu n=\partial_\nu c=0,&x\in\partial\Omega,\ t>0,\\
		(n,c)(x,0)=(n_0,c_0),&x\in\Omega.
	\end{cases}
\end{equation}
Here \(n\) and \(c\) denote the population density and signal concentration, respectively, and \(\nu\) is the outward unit normal vector of the boundary $\partial\Omega$. This problem has a nice chemotactic free energy
\begin{equation}\label{eq:energy}
	\mathscr{E}(n,c):=\int_\Om n\ln n-\int_\Om nc
	+\frac12\int_\Om\left(|\nabla c|^2+c^2\right),
\end{equation}
which enjoys 
\begin{align}\label{eq:0037-1}
	\frac{d}{dt}\mathscr{E}(n,c)+\into n|\nabla(\ln n-c)|^2+\into|\Delta c-c+n|^2=0.
\end{align}
Here and in what follows, whenever no ambiguity arises, we write
\(\int_{\Omega} f\) for \(\int_{\Omega} f(x)\,dx\).

{In light of the observation \eqref{eq:0037-1}}, the celebrated theorem of Nagai, Senba, and Yoshida~\cite{NagaiSenbaYoshida1997} gave global classical solvability of problem \eqref{eq:PP-Neumann} whenever
\begin{equation}\label{eq:PP-classical-thresholds}
	m<\begin{cases}
		4\pi,&\text{on a bounded planar domain with smooth boundary},\\
		8\pi,&\text{for radial data on a disk}.
	\end{cases}
\end{equation}
Its enduring strength lies in converting free-energy dissipation and the optimal Trudinger-Moser inequality into an explicit mass condition. Sharpness of the  threshold \eqref{eq:PP-classical-thresholds} was established by the influential nonsymmetric analysis of Horstmann and Wang~\cite{HorstmannWang2001}: on simply connected planar domains, suitable data of every mass \(m\in(4\pi,\infty)\setminus4\pi\mathbb{N}\) generate unbounded solutions. The radial counterpart on a disk exhibits the same mass-sensitive structure: the groundbreaking collapse analysis~\cite{HerreroVelazquez1997} identifies finite-time concentration for initial mass \(m>8\pi\).

The corresponding mass scale for the Cauchy problem in \(\mathbb{R}^2\) is \(8\pi\). Global existence and boundedness for \(m<8\pi\) are established in~\cite{Mizoguchi2013} and \cite{Hosono2026Radial}, respectively. A remarkable endpoint theorem in~\cite{Hosono2026} extends global existence to \(m=8\pi\). A recent substantial advance~\cite{Hosono2026Radial} addresses precisely this distinction. If the equation of signal is replaced by $c_t=\Delta c+n$, radial blow-up constructions and stable blow-up dynamics in the supercritical range \(m>8\pi\) are developed in~\cite{Mizoguchi2020Cauchy,Schweyer2014}.

Despite these major advances, the corresponding results for the parabolic-parabolic Keller-Segel system \eqref{eq:PP-Neumann} remain less extensive than those for its parabolic-elliptic counterpart. The latter has two important variants: the second equation in \eqref{eq:PP-Neumann} is replaced by the elliptic equation \(0=\Delta c-c+n\) in bounded domains or by \(0=\Delta c+n\) in the whole plane. Although the critical mass of such a parabolic-elliptic Keller-Segel system coincides with that of the fully parabolic case, the fact that the signal is governed by an elliptic equation permits the application of more specialized mathematical techniques, thereby yielding a broader range of results (see \cite{Nagai2001,BlanchetDolbeaultPerthame2006,BlanchetCarrilloMasmoudi2008,BlanchetCarlenCarrillo2012,Wei2018,CarlenFigalli2013,GhoulMasmoudi2018,DavilaDelPinoDolbeaultMussoWei2024,RaphaelSchweyer2014,CollotGhoulMasmoudiNguyen2022Stable,CollotGhoulMasmoudiNguyen2022Spectral,Mizoguchi2022Refined,BuseghinDavilaDelPinoMusso2026} and references therein).

Let us point out that when the signal satisfies an elliptic equation, it is determined by the density via Green-function representations, thereby reducing the system to a single nonlocal evolution equation. By contrast, when the signal satisfies a parabolic equation, it retains its independent initial datum and evolves dynamically in response to the density, so estimates based on elliptic elimination do not transfer directly. Moreover, as highlighted by Winkler~\cite{Winkler2020SinglePointBlowUp}, the directional effects of the cross-diffusion in Keller-Segel systems give rise to substantially more intricate mechanisms than those arising in typical scalar reaction-diffusion equations, calling for more delicate analytical arguments. Consequently, detecting the critical mass in the fully parabolic setting is already a substantial challenge even before fluid feedback is introduced.
 
\subsection{Fluid effects on the critical mass in the Keller-Segel system}\label{sec:parabolic-fluid}

Fluid transport introduces a further mechanism into the competition between diffusion and chemotactic attraction. The influential perspectives~\cite{BellomoTaoWinkler2020,BellomoOutadaSolerTaoWinkler2022} place this interaction within the broader theory of chemotaxis in complex environments.

It is shown that when the velocity field is prescribed as a passive advection—thus omitting dynamic feedback from the population—the strategic design of the flow can suppress singularity formation and raise the effective critical mass. The pioneering mixing results of~\cite{KiselevXu2016} and the powerful shear-flow analyses of~\cite{BedrossianHe2017} establish that suitably chosen incompressible flows, such as relaxation-enhancing flows or shear flows~\cite{ConstantinKiselevRyzhikZlatos2008,YaoZlatos2017,MasmoudiZhao2024}, can suppress blow-up at arbitrary mass in two-dimensional periodic settings; on the plane, an ingenious fast-splitting mechanism~\cite{HeTadmor2019} yields global evolution for mass \(m<16\pi\).

Recent studies demonstrate that blow-up can also be suppressed when the fluid flow is dynamically coupled to the population density through a population-dependent forcing term. For Navier-Stokes coupling, a notable advance~\cite{ZengZhangZi2021} established global solvability at arbitrary mass near a sufficiently strong Couette flow. A corresponding stabilization result near Poiseuille flow is obtained in~\cite{LiXiangXu2025} under an additional smallness condition on the initial vorticity perturbation. These results rely primarily on the enhanced dissipation mechanism in the hydrodynamic stability studies \cite{MasmoudiZhao2024,WeiZhang2026}. Other feedback couplings, such as those in~\cite{HuKiselev2024,HuKiselevYao2025}, reveal further striking regularization mechanisms, underscoring the importance of the specific fluid law and forcing mechanism.

A decisive advance in the critical mass theory for active fluid coupling is the conclusion in~\cite{WangWeiZhang2026} for
\begin{equation*}
	\begin{cases}
		n_t+u\cdot\nabla n=\Delta n-\nabla\!\cdot(n\nabla c), & x\in\mathbb{R}^2,\;t>0,\\
		-\Delta c=n, & x\in\mathbb{R}^2,\;t>0,\\
		u_t+(u\cdot\nabla)u=\Delta u-\nabla P+n\nabla\phi, & x\in\mathbb{R}^2,\;t>0,\\
		\nabla\!\cdot u=0, & x\in\mathbb{R}^2,\;t>0.
	\end{cases}
\end{equation*}
Under suitable regularity assumptions on the initial data and boundedness of \(\nabla\phi\), global strong solvability is established whenever $\int_{\mathbb{R}^2}n_0\le8\pi$. Building on the Green-function representation of the density-determined signal, a key ingredient of their approach is the rearrangement of the density function $n$ and the derivation of its governing differential inequality. The remarkable strength of the rearrangement argument lies in obtaining density estimates independently of the velocity, thereby preserving sharp mass information despite nonlinear fluid feedback. This result shows that the critical mass for global existence and blow-up in the Keller-Segel-Navier-Stokes system on the whole plane  \(\mathbb{R}^2\) is again \(8\pi\), without requiring a dominant shear or a prescribed stabilizing flow geometry.

Nevertheless, sharp mass-dependent results for dynamically evolving signal components remain comparatively scarce in contrast to the setting where the signal satisfies an elliptic equation. This is especially true for fluid flows without a prescribed dominant flow structure. The underlying difficulty is structural: the signal evolves in time and is advected by the fluid, which invalidates the rearrangement arguments for the unknown density that have been exploited in prior works. Moreover, interactions between chemical gradients and velocity gradients break the monotonicity of the classical chemotactic free energy. On bounded domains, potential boundary concentration poses an additional obstacle to retaining the sharp mass threshold. As a result, the existing theory leaves a substantial open gap concerning whether the critical mass for the fully parabolic two-dimensional Keller-Segel system coupled to the Navier-Stokes equations remains $4\pi$, as in its fluid-free counterpart.

\subsection{Main result}\label{sec:main-result}

Motivated by the above works, we consider the  initial-boundary value problem
\begin{equation}\label{eq:KSNS}
	\left\{
	\begin{aligned}
		&n_t+u\cdot\nabla n=\Delta n-\nabla\!\cdot(n\nabla c),
		&&x\in\Om,\ t>0,\\
		&c_t+u\cdot\nabla c=\Delta c-c+n,
		&&x\in\Om,\ t>0,\\
		&u_t+(u\cdot\nabla)u=\Delta u-\nabla P+n\nabla\phi,
		&&x\in\Om,\ t>0,\\
		&\nabla\!\cdot u=0,
		&&x\in\Om,\ t>0,\\
		&\partial_\nu n=\partial_\nu c=0,\quad u=0,
		&&x\in\partial\Om,\ t>0,\\
		&n(x,0)=n_0(x),\quad c(x,0)=c_0(x),\quad u(x,0)=u_0(x),
		&&x\in\Omega
	\end{aligned}
	\right.
\end{equation}
in a bounded domain $\Om\subset\R^2$ with smooth boundary. Here $u$ and $P$ are the fluid velocity and pressure, respectively, $\phi=\phi(x)$ is a prescribed smooth potential, and $\nu$ is the outward normal vector of the boundary $\partial\Omega$. Both cells and signal are transported by the fluid, while $n\nabla\phi$ describes the feedback through buoyancy. For the study of global solvability in such a coupled system, the central structural obstacle is the loss of classical Lyapunov identity. Compared with \eqref{eq:0037-1}, the free-energy defined as in \eqref{eq:energy} satisfies
\begin{equation}\label{eq:energydefect}
	\frac{d}{dt}\mathscr{E}(n,c)
	+\int_\Om n\bigl|\nabla(\ln n-c)\bigr|^2
	+\int_\Om|\Delta c-c+n|^2
	=-\int_\Om \nabla c\cdot(\nabla u\cdot\nabla c).
\end{equation}
The right-hand-side term has no sign and couples the chemical gradient directly to the velocity gradient. The difficulty posed by this strongly coupled nonlinear term cannot be resolved merely by matching the Navier-Stokes kinetic energy, thereby forcing us to seek new (quasi-)energy structures to address it.

In a landmark contribution, Winkler~\cite{Winkler2020SIAM} established the first explicit mass range $m<2\pi$ for this fully coupled problem \eqref{eq:KSNS} and he constructed a global generalized solution without any smallness assumptions on the initial signal or velocity. Moreover, for the narrower range \(m<m_*\), where  \(m_*\in(0,2\pi)\), the solution eventually becomes classical and converges to \((\bar{n}_0,\bar{n}_0,0)\) as $t\to\infty$. Here $\overline{n}_0:=\frac{1}{|\Omega|}\into n_0$. This work decisively established mass as an effective control quantity despite~\eqref{eq:energydefect}. Even so, the remaining gap is equally precise: the explicit mass bound for global solvability is only half of the corresponding bound for the classical Keller-Segel system, and the solution concept does not provide classical regularity at all positive times.  {As is pointed out by Winkler in \cite[p. 2045]{Winkler2020SIAM},  this result leaves open the possibility that fluid interaction in the Keller-Segel-Navier-Stokes system might sightly alter the critical blow-up mass of the classical Keller-Segel system, potentially reducing it to a level between $2\pi$ and $4\pi$.}


{The main purpose of this work is to resolve Winkler's open problem regarding the sharpness of the mass threshold for the initial-boundary value problem \eqref{eq:KSNS}. Specifically, {\bf we establish uniform-in-time boundedness for the system whenever the total mass \(m<4\pi\), without requiring any smallness assumptions on the initial signal or fluid velocity.}}

To state the result precisely, we first introduce the functional setting and initial-data assumptions. Let $L^2_\sigma(\Om)$ be the solenoidal subspace of $L^2(\Om;\R^2)$ and let $\mathcal{A}=-\mathcal P\Delta$ be the Stokes operator, where $\mathcal P$ is the Helmholtz projection, with $D(\mathcal{A})=W^{2,2}(\Om;\R^2)\cap W^{1,2}_0(\Om;\R^2)\cap L^2_\sigma(\Om)$ (\cite[Section III.2]{Sohr2001}), and for $\alpha>0$ we let $\mathcal{A}^\alpha$ represent the corresponding sectorial fractional powers (\cite[Section 1.4]{Henry1981}). We assume
\begin{equation}\label{eqn-conditionofphi}
 \phi\in W^{2,\infty}(\Om)
\end{equation}
and
\begin{equation}\label{eqn-conditionofinitialda}
 \begin{cases}
 n_0\in C^0(\overline\Om),\quad n_0\ge0,\quad n_0\not\equiv0,\\
 c_0\in W^{1,\infty}(\Om),\quad c_0\ge0,\\
 u_0\in D(\mathcal{A}).
 \end{cases}
\end{equation}
We are now in a position to state our main result.
\begin{theorem}\label{thm:main}
 Let $\Om\subset\R^2$ be a bounded domain with smooth boundary, and assume~\eqref{eqn-conditionofphi}-\eqref{eqn-conditionofinitialda} hold. If
 \begin{equation}\label{eq:newthreshold}
 m:=\int_\Om n_0<4\pi,
 \end{equation}
 then there is a unique triple $(n,c,u)$ enjoying
 \begin{equation*}
 \left\{\begin{aligned}
 n&\in C^0(\overline\Om\times[0,\infty))\cap C^{2,1}(\overline\Om\times(0,\infty)),\\
 c&\in\bigcap_{q>2}C^0([0,\infty);W^{1,q}(\Om))\cap C^{2,1}(\overline\Om\times(0,\infty)),\\
 u&\in\bigcap_{\alpha\in(\tfrac{1}{2},1)}C^0([0,\infty);D(\mathcal{A}^\alpha))
 \cap C^{2,1}(\overline\Om\times(0,\infty);\R^2),
 \end{aligned}\right.
 \end{equation*}
 such that~\eqref{eq:KSNS} is satisfied classically with some $P\in C^{1,0}(\Om\times(0,\infty))$. Moreover, this solution is uniformly bounded in time: for every $q>2$, there is $C>0$ such that
\begin{equation}\label{eq:uniform-bounds}
	\|n(\cdot,t)\|_{L^\infty(\Omega)}+\|c(\cdot,t)\|_{W^{1,q}(\Omega)}+\|u(\cdot,t)\|_{L^\infty(\Omega)}\le C\quad\text{for all }t>0.
\end{equation}
\end{theorem}

\begin{remark}\label{rem:gap}
	Compared with the result obtained in \cite{Winkler2020SIAM}, the improvement is twofold: the mass restriction is relaxed from \(2\pi\) to \(4\pi\), and the conclusion is upgraded {from generalized solvability} and eventual regularity under an additional smallness condition {to global classical solvability} for all positive times with uniform-in-time boundedness. In particular, the theorem rules out the reduction of the critical mass to a level between \(2\pi\) and \(4\pi\) contemplated in the discussion recalled above. {It therefore resolves Winkler's open problem and affirmatively answers the conjecture of Bonilla and Guti\'errez-Santacreu~\cite{BonillaGutierrezSantacreu2024} that the classical subcritical range is preserved under fluid coupling.}
\end{remark}

\begin{remark}\label{rem:sharpness}
	It should be noted that the threshold \(4\pi\) in Theorem~\ref{thm:main} is optimal in the sense that the conclusion fails for suitable data with \(m>4\pi\). Indeed, in a simply connected planar domain with smooth boundary, take \(\phi\) to be constant and \(u_0=0\). Then \(u\equiv0\), and~\eqref{eq:KSNS} reduces to~\eqref{eq:PP-Neumann}. For each \(m\in(4\pi,\infty)\setminus4\pi\mathbb N\), Horstmann and Wang~\cite{HorstmannWang2001} constructed data for~\eqref{eq:PP-Neumann} whose density becomes unbounded either in finite time or as \(t\to\infty\).
\end{remark}

\begin{remark}
	To the best of our knowledge, the critical mass case for the problem \eqref{eq:KSNS} remains still open, even for the fluid-free counterpart. The reason is that the coercivity in Lemma \ref{lem:lowerboundofentropy} degenerates at \(m=4\pi\) and the coefficient in \eqref{eq:delta} vanishes in this case. The critical case requires a different approach and new estimates, and will be addressed in our forthcoming work.
\end{remark}

\subsection{Proof strategy and organization of the paper}\label{sec:proof-idea}

It is well known that the energy method provides a standard and effective approach for establishing global existence results in evolution equations. As noted above, unlike the Keller-Segel system \eqref{eq:PP-Neumann}, the classical Lyapunov identity is no longer available in the present coupled Keller-Segel-fluid setting, and the free energy $\mathscr E(t)$ defined in \eqref{eq:energy} actually satisfies the following identity (Lemma \ref{lem:evolutionofe})
\begin{equation}\label{eq:energydefect-1}
	\frac{d}{dt}\mathscr{E}(t)
	+\int_\Om n\bigl|\nabla(\ln n-c)\bigr|^2
	+\int_\Om|\Delta c-c+n|^2
	=-\int_\Om \nabla c\cdot(\nabla u\cdot\nabla c).
\end{equation}
Our main task is to control the right-hand-side term so as to close the energy estimates.

The novelty of our proof lies, on the one hand, in an interpolation inequality with explicit leading coefficients and, on the other hand, in the key observation of a quasi-energy structure suitable for the system. In particular, the relative isoperimetric geometric properties and rearrangement techniques play a crucial role in establishing this key interpolation inequality.

 Noticing that one can deduce from the fluid equation that (Lemma \ref{lem:2104-2})
\begin{align}\label{eqn:2104-2-1}
	\frac{d}{dt}\int_\Omega |u|^2+\int_\Omega|\nabla u|^2\le C_1\ln\left(\int_\Omega n^2+1\right)+C_1
\end{align}
for some constant $C_1>0$ independent of $t$, we devote our efforts to deriving the following estimate (Lemma \ref{lem:static}):
\begin{equation}\label{eq:idea-static}
	\|n\|_{L^2(\Omega)}^2+\|c\|_{W^{2,2}(\Omega)}^2
	\le C_2\big(1+\mathscr D\big)
\end{equation}
with some $C_2>0$ independent of $t$, which holds under the mass threshold \(m<4\pi\). Here
\begin{equation*}
	\mathscr D:=\int_\Omega n|\nabla(\ln n-c)|^2+\int_\Omega|\Delta c-c+n|^2
\end{equation*}
denotes the dissipation appearing in \eqref{eq:energydefect-1}.

Substituting \eqref{eq:idea-static} into \eqref{eqn:2104-2-1} yields some $C_3>0$ such that
\begin{align*}
	\frac{d}{dt}\int_\Omega|u|^2+\int_\Omega|\nabla u|^2\le C_3\ln(1+\mathscr{D})+C_3.
\end{align*}
Introducing the functional $\mathscr{F}=\log(\mathcal M+\mathscr E)+(\ell+1)\|u\|_{L^2(\Om)}^2$ for some suitable $\mathcal{M},\ell>0$ , we can then derive a new quasi-energy inequality:
\begin{equation*}
	\mathscr{F}'+\frac14\cdot\frac{\mathscr D}{\mathcal M+\mathscr E}+\|\nabla u\|_{L^2(\Om)}^2\le C_4\mathscr{F},
\end{equation*}
where $C_4>0$ is independent of $t$ (Lemma \ref{lem:evolutionofcoupledentroenergy}). From this quasi-energy inequality we obtain preliminary bounds for \(\int_\Omega n\ln n\) and \(\int_\Omega|\nabla c|^2\) (Lemma~\ref{lem:boundofllogl}). Combined with standard bootstrapping arguments and semigroup techniques, these bounds eventually lead to the global existence and boundedness of solutions (Propositions~\ref{lem:boundednessforu}-\ref{lem:boundednessforncu}).

It is worth emphasizing that the proof of \eqref{eq:idea-static} is far from obvious. The construction relies crucially on a spatial interpolation estimate: for every $\varepsilon>0$, there is $C_5>0$ such that
\begin{equation}\label{eq:idea-GNS}
	\int_\Om f^6\le\left(\frac2\pi+\varepsilon\right)
	\left(\int_\Om f^4\right)\left(\int_\Om|\nabla f|^2\right)
	+C_{5}\left(\int_\Om f^4\right)^{\frac{3}{2}}
\end{equation}
for all nonnegative \(f\in C^1(\overline\Om)\) (Lemma~\ref{lem:GNS}). Inspired by \cite{BlanchetCarlenCarrillo2012}, we apply this new interpolation inequality to the function \(f_k=(n+k)^{\frac14}-k^{\frac14}\) with an adaptive shift parameter \(k=\eta\int_\Omega n^2\) for some sufficiently small \(\eta>0\) (Lemmas \ref{lem:staticabsorbstart}-\ref{lem:static}) to establish \eqref{eq:idea-static}. The interpolation inequality \eqref{eq:idea-GNS} is obtained by transporting the sharp planar Gagliardo-Nirenberg-Sobolev inequality \cite{DPD02,CarlenFigalli2013} to bounded domains via symmetrization/rearrangement, as detailed in Lemmas \ref{lem:propertyofrearragement}, \ref{lem:GNS}. This step rests on a relative isoperimetric estimate
$$\mathscr P_\Om(E)^2\ge(2\pi-\varepsilon)|E|,$$
valid for sets $E$ of sufficiently small measure (Lemma~\ref{lem:reliso}), where $\mathscr P_\Om(E)$ denotes the relative perimeter of $E$ in $\Omega$ and $|E|$ represents the Lebesgue measure of $E$.

Isoperimetric inequalities constitute a powerful source of sharp estimates for partial differential equations. Their broad applicability is demonstrated by striking eigenvalue optimization results for drift-diffusion elliptic operators in \cite{HamelNadirashviliRuss2011} and sharp Robin comparison principles in \cite{AlvinoNitschTrombetti2023}, which convert geometric information into quantitative spectral or solution bounds. While sharp functional inequalities and rearrangement methods are well established tools in chemotaxis theory \cite{CarlenFigalli2013,WangWeiZhang2026}, the direct application of small-volume relative isoperimetry to Keller-Segel-fluid systems remains largely unexplored. In our proof, this geometric input preserves the leading constant in the interpolation bound and thereby unlocks the full admissible range \(m<4\pi\).

{\bf Organization of the paper.} {Section~\ref{sec:preliminaries} first recalls the local existence theory and the continuation criterion. It then develops the relative isoperimetric and interpolation estimates. Several useful results are recalled in this section.} Section~\ref{sec:global} implements these tools in the coupled energy analysis: we derive the dissipation estimate and establish the quasi-energy structure. Higher-order regularity estimates and the complete proof of Theorem~\ref{thm:main} are given in Section~\ref{sec:continuation}. Appendix~\ref{sec:rearrangement} collects the rearrangement properties used in the interpolation argument.

\section{Preliminaries}\label{sec:preliminaries}

{This section reviews the local existence theory and several useful analytical tools. In addition, we derive relative isoperimetric and interpolation estimates.}

\subsection{Local existence and extensibility}

The following basic result on the local existence of classical solutions, as well as a convenient continuation criterion, can be derived by adapting standard arguments (cf. \cite{Winkler2012cpde,W19nodea,WangWinkler2024} for closely related settings).
\begin{lem}\label{lem:localexistence}
	Assume that \eqref{eqn-conditionofphi} and \eqref{eqn-conditionofinitialda} hold. Then there exist $\Tmax\in(0,\infty]$ and a unique triple $(n,c,u)$ of functions
	\begin{align*}
		\left\{\begin{aligned}
			n&\in C^0(\overline\Om\times[0,\Tmax))
			\cap C^{2,1}(\overline\Om\times(0,\Tmax)),\\
			c&\in\bigcap_{q>2}C^0([0,\Tmax);W^{1,q}(\Om))
			\cap C^{2,1}(\overline\Om\times(0,\Tmax)),\quad\text{and}\\
			u&\in\bigcap_{\frac{1}{2}<\alpha<1}C^0([0,\Tmax);D(\mathcal{A}^\alpha))
			\cap C^{2,1}(\overline\Om\times(0,\Tmax);\R^2),
		\end{aligned}\right.
	\end{align*}
	such that the initial-boundary value problem \eqref{eq:KSNS} holds in $\Omega\times(0,T_{\max})$ with some $P\in C^{1,0}(\Om\times(0,\Tmax))$. Moreover, we have that
	\begin{align}\label{eqn:extensibility}
		&\text{if }T_{\max}<\infty,\;\text{then}\;\lim_{t\nearrow T_{\max}}\left(\|n(\cdot,t)\|_{L^\infty(\Omega)}+\|c(\cdot,t)\|_{W^{1,q}(\Omega)}+\|\mathcal{A}^\alpha u(\cdot,t)\|_{L^2(\Omega)}\right)=\infty
	\end{align}
	for $q>2$ and $\alpha\in(\frac{1}{2},1)$.
\end{lem}

Mass conservation for the density follows directly from integrating the density equation, whereas the \(L^1\) estimate for \(c\) and the positivity of \(n\) and \(c\) are obtained from the comparison principle.
\begin{lem}
	We have that $n>0$ and $c>0$ in $\overline{\Omega}\times(0,T_{\max})$. Moreover,
	\begin{align}\label{eqn-massconservation}
		\into n(\cdot,t)=\into n_0\quad\text{for all }t\in(0,T_{\max}),
	\end{align}
	and
	\begin{align}\label{eqn-signalmass}
		\into c(\cdot,t)\le\max\left\{\into n_0,\into c_0\right\}\quad\text{for all }t\in(0,T_{\max}).
	\end{align}
\end{lem}

\subsection{An interpolation inequality resulting from relative isoperimetry}

To establish global existence for the system, we derive an interpolation inequality with an explicit leading coefficient in this subsection. This will be accomplished using a relative isoperimetric estimate and rearrangement arguments.

\subsubsection{A relative isoperimetric estimate at small volumes}

The relevant geometric quantity is the perimeter measured inside the domain. For a Lebesgue measurable set $E \subset \Omega$, its \emph{relative perimeter} in $\Omega$ is
\[
	\mathscr{P}_\Omega(E):= \sup\left\{\int_E\operatorname{div} X\,dx: X\in C_0^1(\Omega;\R^2),|X|\le 1\text{ in }\Omega\right\},
\]
which is the standard variational definition of the perimeter of $E$ relative to the open set $\Omega$ (see \cite[Chapter 12]{Maggi2012} or \cite[Definition 5.1.2]{BFR26} for instance).

Following \cite[\S 2]{BayleRosales2005} (see also \cite[\S 2.4]{Ros2005}), we define the \emph{(relative) isoperimetric profile} of $\Omega$ as
\begin{equation}\label{def of isoperimetric profile}
	\mathscr{I}_\Omega(l):=\inf\left\{\mathscr{P}_\Omega(E):E \subset\Omega \text{ measurable},|E|=l\right\},\quad 0<l<|\Omega|.
\end{equation}
We say that $E$ is a set of finite perimeter in $\Omega$ if $\mathscr{P}_\Omega(E)<\infty$.
\begin{lem}\label{lem:reliso}
	Let $\Omega \subset \mathbb{R}^2$ be a bounded domain with smooth boundary. For every $\varepsilon > 0$ there exists $\mathfrak{s}\in (0, |\Omega|]$ such that
	\begin{equation}\label{eq:reliso}
		\mathscr{P}_\Omega(E)^2 \ge (2\pi - \varepsilon)\,|E|
	\end{equation}
	for every finite-perimeter set $E \subset \Omega$ with $|E| \le \mathfrak{s}$.
\end{lem}
\begin{proof}
	Thanks to Theorem 1.3 in \cite{Ros2005} (see also \cite[\S 2]{BayleRosales2005}), the isoperimetric profile of the half-space $\mathbb{R}^{d}_+:=\{x=(x_1,x_2,\cdots,x_d)\in\mathbb{R}^d:\,x_{d}>0\}$ is
	\[
		\mathscr{I}_{\mathbb{R}_+^d}(l)= 2^{-\frac{1}{d}}\gamma_{d}\,l^{\frac{d-1}{d}},
	\]
	where $l>0$ and $\gamma_{d}:= \frac{\mathcal{H}^{d-1}(\mathbb{S}^{d-1})}{|B_1|^{\frac{d-1}{d}}}$. The notation $\mathcal{H}^k$ denotes the $k$-dimensional Hausdorff measure and $\mathbb{S}^k$ is the unit sphere centered at the origin in $\mathbb{R}^{k+1}$. In particular, for $d = 2$, we have $\gamma_2 = 2\sqrt{\pi}$ and hence
	\begin{equation*}
		\mathscr{I}_{\mathbb{R}^2_+}(l) = \sqrt{2\pi l},\quad l>0.
	\end{equation*}
	The small-volume asymptotics in \cite[Proposition~2.1 and Remark~2.2]{BayleRosales2005} or \cite[Corollary~1.3]{Fall2010} therefore show that
	\begin{equation*}
		\mathscr{I}_\Omega(\ell) \sim \mathscr{I}_{\mathbb{R}^2_+}(l) = \sqrt{2\pi l}
		\quad \text{as }l\searrow 0.
	\end{equation*}
	Consequently, for every $\eta>0$, there exists $\widetilde{l}>0$ such that
	\begin{equation}\label{eq:twosided}
		(1-\eta)\sqrt{2\pi l}\le \mathscr{I}_\Omega(l) \le(1+\eta)\sqrt{2\pi l}\quad \text{for all } 0 < l \le \widetilde{l}.
	\end{equation}
	
	It suffices to consider $0<\varepsilon<2\pi$, since the conclusion is trivial for $\varepsilon\ge2\pi$. Choose $\eta := \varepsilon/(4\pi)$. For every
	finite-perimeter set $E \subset \Omega$ with
	$l:=|E|\le\widetilde{l}$, according to the definition of the profile \eqref{def of isoperimetric profile} and the lower bound in \eqref{eq:twosided}, we know
	\[
		\mathscr{P}_\Omega(E) \ge \mathscr{I}_\Omega(l) \ge (1-\eta)\sqrt{2\pi l},
	\]
	which, upon squaring, yields
	\[
		\mathscr{P}_\Omega(E)^2\ge(1-\eta)^22\pi|E|\ge (1-2\eta)2\pi|E|=(2\pi-\varepsilon)|E|.
	\]
	Hence \eqref{eq:reliso} holds with $\mathfrak{s}:= \min\{ \widetilde{l}, |\Omega|\}$.
\end{proof}

\subsubsection{Rearrangement}

Let $\Omega\subset\mathbb{R}^2$ be a bounded domain with smooth boundary. For any $a\in(0,2\pi)$, we use Lemma \ref{lem:reliso} to fix {$\mathfrak{s}\in(0,|\Omega|)$, decreasing it if necessary}. For this choice of $\mathfrak{s}$ and a nonnegative function $0\not\equiv f\in C^1(\overline{\Omega})$, define
\begin{align}\label{eqn:def-1}
	b:=\left(\frac{2}{\mathfrak{s}}\into f^4\right)^{\frac{1}{4}},\quad g(x):=(f(x)-b)_+,\quad S:=|\{x\in\Omega:g(x)>0\}|.
\end{align}
This choice shows $S\le \frac{\mathfrak{s}}{2}$, while the gradient of the truncated function remains controlled by that of $f$. We also define the decreasing rearrangement\footnote{The distribution-function identities, equimeasurability and regularity properties needed here are recalled in Appendix~\ref{sec:rearrangement}.} by
\begin{align}\label{eqn:def-2}
	g^*(s):=\inf\{r\ge 0:|\{x\in\Omega:g(x)>r\}|\le s\},\quad \text{for }s\in(0,|\Omega|).
\end{align}
The following estimate converts the relative perimeter bound into control of the weighted derivative of $g^*$.

\begin{lem}\label{lem:propertyofrearragement}
	Let $\Omega\subset\mathbb{R}^2$ be a bounded domain with smooth boundary. For any $a\in(0,2\pi)$, one can find $\mathfrak{s}\in(0,|\Omega|)$ such that the quantities $g,S,g^*$, defined by \eqref{eqn:def-1}-\eqref{eqn:def-2}, satisfy
		\begin{align}\label{eq:1558-1}
			\int_0^Ss\left|(g^*)'(s)\right|^2ds\le\frac{1}{a}\into|\nabla g|^2.
		\end{align}
\end{lem}
\begin{proof}
	For any fixed $a\in(0,2\pi)$, Lemma \ref{lem:reliso} provides $\mathfrak{s}$ such that
	\begin{align}\label{eqn-1548-1}
		\mathscr{P}_\Omega(E)^2\ge a|E|\quad\text{for each }E\subset\Omega\;\text{of finite relative perimeter fulfilling }0<|E|\le \mathfrak{s}.
	\end{align}
	For $0<s<t\le S$, we infer from \eqref{eqn:1451-2} and \eqref{eqn:1522-1} in Lemma \ref{lem:some basic properties of the rearrangement} that
	\begin{align*}
		s<|\{x\in\Omega:g>r\}|\le\mathfrak{s}\quad\text{for }r\in(g^*(t),g^*(s)),
	\end{align*}
	and
	\begin{align*}
		\int_{g^*(t)}^{g^*(s)}\mathscr{P}_\Omega(\{x\in\Omega:g>r\})dr=\int_{\{x\in\Omega:g^*(t)<g<g^*(s)\}}|\nabla g|.
	\end{align*}
	Applying \eqref{eqn-1548-1} to $\{x\in\Omega:g>r\}$, we arrive at
	\begin{align*}
		(as)^{\frac{1}{2}}\left(g^*(s)-g^*(t)\right)&\le \int_{g^*(t)}^{g^*(s)}\mathscr{P}_\Omega(\{x\in\Omega:g>r\})dr\nonumber\\
		&=\int_{\{x\in\Omega:g^*(t)<g<g^*(s)\}}|\nabla g|.
	\end{align*}
	Since \eqref{eqn:1451-1} in Lemma \ref{lem:some basic properties of the rearrangement} shows $|\{x\in\Omega:g^*(t)<g<g^*(s)\}|\le t-s$, we apply H\"older's inequality to see that
	\begin{align*}
		(as)^{\frac{1}{2}}\left(g^*(s)-g^*(t)\right)\le(t-s)^{\frac{1}{2}}\left(\int_{\{x\in\Omega:g^*(t)<g<g^*(s)\}}|\nabla g|^2\right)^{\frac{1}{2}}.
	\end{align*}
	Setting $t:=s+h$ with $h>0$ sufficiently small and integrating the resulting difference-quotient estimate, we obtain
	\begin{align*}
		&a\int_0^{S-h}s\left|\frac{g^*(s+h)-g^*(s)}{h}\right|^2ds\\
		&\quad\le\frac{1}{h}\left(\int_0^{S-h}\int_{\{x\in\Omega:g>g^*(s+h)\}}|\nabla g|^2dxds-\int_0^{S-h}\int_{\{x\in\Omega:g>g^*(s)\}}|\nabla g|^2dxds\right)\\
		&\quad=\frac{1}{h}\left(\int_h^S\int_{\{x\in\Omega:g>g^*(s)\}}|\nabla g|^2dxds-\int_0^{S-h}\int_{\{x\in\Omega:g>g^*(s)\}}|\nabla g|^2dxds\right)\\
		&\quad=\frac{1}{h}\left(\int_{S-h}^S\int_{\{x\in\Omega:g>g^*(s)\}}|\nabla g|^2dxds-\int_{0}^{h}\int_{\{x\in\Omega:g>g^*(s)\}}|\nabla g|^2dxds\right)\\
		&\quad\le\into|\nabla g|^2,
	\end{align*}
	which, combined with Fatou's lemma, implies \eqref{eq:1558-1} as $h\searrow0$.
\end{proof}

\subsubsection{Interpolation}

Using Lemma \ref{lem:propertyofrearragement}, we are now able to localize the sharp planar Gagliardo-Nirenberg-Sobolev inequality~\cite{DPD02,CarlenFigalli2013} to the bounded domain by means of the Schwarz symmetrization technique.

\begin{lem}\label{lem:GNS}
	Let $\Omega\subset\mathbb{R}^2$ be a bounded domain with smooth boundary. For every $\varepsilon>0$ there is $C>0$ such that
	\begin{align}\label{eq:GNSdomain}
		\int_\Omega f^6\le\left(\frac{2}{\pi}+\varepsilon\right)\cdot\left(\int_\Omega f^4\right)\cdot\left(\int_\Omega|\nabla f|^2\right)+C\cdot\left(\int_\Omega f^4\right)^{\frac{3}{2}}
	\end{align}
	for all nonnegative $f\in C^1(\overline{\Omega})$.
\end{lem}
\begin{proof}
	We first recall the sharp planar Gagliardo-Nirenberg-Sobolev inequality (cf. \cite[Theorem 1]{DPD02} or \cite[equation (1.8)]{CarlenFigalli2013}):
	\begin{align}\label{eq-gns}
		\pi\int_{\mathbb R^2}v^6\le\left(\int_{\mathbb R^2}v^4\right)\cdot\left(\int_{\mathbb R^2}|\nabla v|^2\right)\quad\text{for }v\in W^{1,2}(\mathbb{R}^2).
	\end{align}
	For any $a\in(0,2\pi)$, we apply Lemma \ref{lem:reliso} to choose $\mathfrak{s}\in(0,|\Omega|)$. For this choice of $\mathfrak{s}$, define $b,g,S,g^*$ as in \eqref{eqn:def-1} and \eqref{eqn:def-2}. We then define the Schwarz symmetrization associated with $g^*$ by
	\[
	G(x):=
	\begin{cases}
		g^*(\pi|x|^2),&\text{if }0<\pi|x|^2<S,\\
		0,&\text{if }\pi|x|^2\ge S{\text{ or }x=0}.
	\end{cases}
	\]
	Applying (iv) in Lemma \ref{lem:some basic properties of the rearrangement} and Lemma \ref{lem:propertyofrearragement}, we see that for any $a\in(0,2\pi)$,
	\begin{equation}\label{eq:radialgradient}
		\int_{\mathbb R^2}|\nabla G|^2\,dx
		=4\pi\int_0^S s|(g^*)'(s)|^2\,ds
		\le\frac{4\pi}{a}\int_\Omega|\nabla g|^2\,dx.
	\end{equation}
	Moreover, equimeasurability \cite[Property 1.1.2 and Corollary 1.1.1]{Rak25} yields, for $p=4,6$,
	\begin{align}\label{eq:rearragementineq}
		\int_{\mathbb R^2}(G(x))^p\,dx=\int_0^S(g^*(s))^p\,ds=\int_\Omega g^p(x)\,dx.
	\end{align}
	Applying \eqref{eq-gns} to $G$ and using \eqref{eq:radialgradient} and \eqref{eq:rearragementineq}, we obtain, for each $a\in(0,2\pi)$,
	\begin{equation}\label{eq:smallGNS}
		\int_\Omega g^6
		\le\frac4a\left(\int_\Omega g^4\right)\cdot
		\left(\int_\Omega|\nabla g|^2\right).
	\end{equation}
	
	Notice that for every $\theta>1$, $s^6\le\theta(s-b)_+^6$ if $s\ge\frac{b}{1-\theta^{-\frac{1}{6}}}$; $s^6\le \frac{1}{(1-\theta^{-\frac{1}{6}})^2}\cdot b^2s^4$ if $0\le s<\frac{b}{1-\theta^{-\frac{1}{6}}}$. Consequently, one can find $C_1>0$ independent of $b$ such that
	\begin{equation}\label{eq:gns32-remove-cutoff}
		f^6\le\theta(f-b)_+^6+C_1b^2f^4.
	\end{equation}
	Recalling that $b:=\left(\frac{2}{\mathfrak{s}}\into f^4\right)^{\frac{1}{4}}$ and combining \eqref{eq:gns32-remove-cutoff} with \eqref{eq:smallGNS}, we obtain, for $a\in(0,2\pi)$ and $\theta>1$,
	\begin{align}\label{eq:gns32-final-before-choice}
		\int_\Omega f^6\le\frac{4\theta}a\left(\int_\Omega f^4\right)\cdot\left(\int_\Omega|\nabla f|^2\right)+C_1\sqrt{\frac{2}{\mathfrak{s}}}\left(\int_\Omega f^4\right)^{\frac{3}{2}}.
	\end{align}
	Finally, for any $\varepsilon>0$, the inequality $\frac{4}{\frac{2}{\pi}+\varepsilon}<2\pi$ allows us to choose $a\in(0,2\pi)$ such that
	\begin{align*}
		\frac{4}{\frac{2}{\pi}+\varepsilon}<a<2\pi,
	\end{align*}
	which implies $\frac{a}{4}\left(\frac{2}{\pi}+\varepsilon\right)>1$. We may therefore choose $\theta$ satisfying
	\begin{align*}
		1<\theta<\frac{a}{4}\left(\frac{2}{\pi}+\varepsilon\right).
	\end{align*}
	Combining this with \eqref{eq:gns32-final-before-choice} leads to \eqref{eq:GNSdomain}.
\end{proof}

\begin{rem}
	A standard extension and mollification argument shows that \eqref{eq:GNSdomain} also holds for every nonnegative $f\in H^1(\Omega)$. We state the lemma for $C^1(\overline{\Omega})$-functions to keep the rearrangement argument separate from the approximation procedure.
\end{rem}


\subsection{Some useful inequalities and regularity estimates}\label{sec:auxiliary-entropy}

{We first recall Winkler's logarithmic-type inequality~\cite[Lemma~3.6]{W19nodea}.}
\begin{lem}\label{lem:logarithmic-entropy}
	Assume $\Omega\subset\mathbb{R}^d$ is a bounded domain with smooth boundary. Let $\varphi\in L^2(\Omega)$ be such that $\varphi\ge 0$ a.e. in $\Omega$. Then
	\begin{align*}
		\into\varphi\ln\varphi\le\left(\into\varphi+1\right)\cdot\ln\left(\into \varphi^2+e\right).
	\end{align*}
\end{lem}

The second inequality comes from \cite[Lemma~2.2]{Winkler2020SIAM}. Its explicit mass dependence will be essential for the coercivity of the chemotactic free energy below $4\pi$.
\begin{lem}\label{lem:entropy-interaction}
	{Suppose that $\Omega\subset\mathbb{R}^2$ is a bounded domain with smooth boundary. For every $\varepsilon>0$ there exists $M=M(\varepsilon,\Omega)>0$ such that, whenever $0\not\equiv\phi\in C^0(\overline\Omega)$ is nonnegative and $\psi\in W^{1,2}(\Omega)$, we have for every $a>0$,}
	\begin{align*}
		\into \phi|\psi|\le\frac{1}{a}\into\phi\ln\frac{\phi}{\overline{\phi}}+\frac{(1+\varepsilon)a}{8\pi}\cdot\left(\into\phi\right)\cdot\left(\into|\nabla\psi|^2\right)+Ma\left(\into\phi\right)\cdot\left(\into|\psi|\right)^2+\frac{M}{a}\into\phi,
	\end{align*}
	where $\overline\phi:=\frac{1}{|\Omega|}\int_\Omega\phi$.
\end{lem}

The following interpolation result is a corollary of \cite[Proposition~1.1]{WangWinkler2024} with $\alpha=0$.
\begin{lem}\label{lem:orlicz-interpolation}
Let $\Omega\subset\mathbb R^2$ be a bounded domain with smooth boundary. Suppose that $\hbar\in C^0([0,\infty))$ is positive and nondecreasing, with $\hbar(s)\to\infty$ as $s\to\infty$. For every $K>0$ and $\eta>0$, there exists $C>0$ such that every $\varphi\in C^1(\overline\Omega)$ satisfying $\varphi>0$ in $\overline\Omega$ and
\[
 \int_\Omega\varphi\hbar(\varphi)\le K
\]
obeys
\[
 \int_\Omega\varphi^{3}
 \le\eta\int_\Omega|\nabla\varphi|^2+C.
\]
\end{lem}

To transfer density control to the fluid velocity, we also recall the following estimate for the forced Navier-Stokes equations from \cite[Lemma~3.4(ii)]{WangWinkler2024}.
\begin{lem}\label{lem:forced-fluid-regularity}
Let $\Omega\subset\mathbb R^2$ be a bounded domain with smooth boundary, and let $\mathcal{A}$ be the Stokes operator on $L^2_\sigma(\Omega)$ with homogeneous Dirichlet boundary conditions. For every $K>0$ and $\beta\in(\frac12,1)$, there exists $C>0$ with the following property. If $T>0$,
\[
\begin{aligned}
 f&\in C^0(\overline\Omega\times[0,T);\mathbb R^2),\\
 v&\in C^0([0,T);D(\mathcal{A}^\beta))
       \cap C^{2,1}(\overline\Omega\times(0,T);\mathbb R^2),\\
 \Pi&\in C^{1,0}(\Omega\times(0,T))
\end{aligned}
\]
satisfy
\[
\begin{cases}
 v_t+(v\cdot\nabla)v=\Delta v-\nabla\Pi+f,
 \quad\nabla\!\cdot v=0,&x\in\Omega,\ 0<t<T,\\
 v=0,&x\in\partial\Omega,\ 0<t<T,
\end{cases}
\]
and
\[
 \|\mathcal{A}^\beta v(\cdot,0)\|_{L^2(\Omega)}^2\le K,
 \qquad
 \sup_{0<t<T}\|f(\cdot,t)\|_{L^2(\Omega)}^2\le K,
\]
then
\[
 \sup_{0<t<T}\|\mathcal{A}^\beta v(\cdot,t)\|_{L^2(\Omega)}^2\le C.
\]
In particular, the constant $C$ is independent of $T$.
\end{lem}

\section{Global classical solvability in the subcritical mass case}\label{sec:global}

{Combining the interpolation and the logarithmic-type inequality from Section~\ref{sec:preliminaries}, we proceed with the energy analysis. The argument begins by tracking the evolution of the chemotactic free energy, after which the density and signal concentration are controlled via dissipation. From the ordinary differential inequality satisfied by the combination of the logarithm of chemotactic free energy and the fluid kinetic energy, we can finally deduce some fundamental regularity estimates for $n$, $c$, and $u$. }
 
\subsection{The chemotactic free energy and the fluid kinetic energy}

The natural starting point is to examine the evolution of classical chemotactic free energy defined as
\begin{align}\label{eqn:def of entropy}
	\mathscr{E}(t):=\into n(\cdot,t)\ln n(\cdot,t)+\frac{1}{2}\into|\nabla c(\cdot,t)|^2+\frac{1}{2}\into c^2(\cdot,t)-\into n(\cdot,t)c(\cdot,t),\quad t>0.
\end{align}
 

 \begin{lem}\label{lem:evolutionofe}
 	We have
 	\begin{align}\label{eqn:entropy}
 		&\frac{d}{dt}\mathscr{E}(t)+\into n(\cdot,t)|\nabla(\ln n(\cdot,t)-c(\cdot,t))|^2+\into|\Delta c(\cdot,t)-c(\cdot,t)+n(\cdot,t)|^2\nonumber\\
 		&\quad=-\int_\Om \nabla c(\cdot,t)\cdot(\nabla u(\cdot,t)\cdot\nabla c(\cdot,t))
 	\end{align}
 	for all $t\in(0,T_{\max})$.
 \end{lem}
 \begin{proof}
 	Multiplying the first equation of \eqref{eq:KSNS} by $\ln n-c$ and integrating by parts, we have
 	\begin{align*}
 		\frac{d}{dt}\into n\ln n-\into n_tc=-\into n|\nabla(\ln n-c)|^2-\into u\cdot\nabla n(\ln n-c)
 	\end{align*}
 	for all $t\in(0,T_{\max})$. Since $\nabla\cdot u=0$ in $\Omega$, $u=0$ on $\partial\Omega$, and $\frac{d}{dt}\into nc=\into n_tc+\into nc_t$, we have
 	\begin{align}\label{eqn:2025-1}
 		\frac{d}{dt}\left(\into n\ln n-\into nc\right)+\into nc_t=-\into n|\nabla(\ln n-c)|^2-\into nu\cdot\nabla c
 	\end{align}
 	for all $t\in(0,T_{\max})$. In view of the second equation of \eqref{eq:KSNS}, one deduces
 	\begin{align*}
 		\frac{1}{2}\frac{d}{dt}\left(\into|\nabla c|^2+\into c^2\right)&=\into\nabla c\cdot\nabla c_t+\into cc_t\nonumber\\
 		&=\into c_t(c-\Delta c)\quad\text{for all }t\in(0,T_{\max}).
 	\end{align*}
 	Combining this identity with \eqref{eqn:2025-1}, substituting the second equation of \eqref{eq:KSNS}, and integrating by parts once more, we obtain
 	\begin{align*}
 		\mathscr{E}'(t)&=-\into c_t(\Delta c-c+n)-\into n|\nabla(\ln n-c)|^2-\into nu\cdot\nabla c\nonumber\\
 		&=-\into|\Delta c-c+n|^2+\into(u\cdot\nabla c)(\Delta c-c+n)-\into n|\nabla(\ln n-c)|^2-\into nu\cdot\nabla c\nonumber\\
 		&=-\into|\Delta c-c+n|^2-\into n|\nabla(\ln n-c)|^2-\int_\Om \nabla c\cdot(\nabla u\cdot\nabla c)
 	\end{align*}
 	for all $t\in(0,T_{\max})$, which is \eqref{eqn:entropy}.
 \end{proof}
 
The right-hand side of \eqref{eqn:entropy} has no prescribed sign. Its control therefore requires information from the fluid equations. For this purpose, the useful feature of the next estimate is the logarithmic dependence of the fluid forcing on $\int_\Omega n^2$: it is compatible with the logarithmic normalization of the free energy used below.

 \begin{lem}\label{lem:2104-2}
 	There is $C>0$ such that
 	\begin{align}\label{eqn:2104-2}
 		\frac{d}{dt}\into |u(\cdot,t)|^2+\into|\nabla u(\cdot,t)|^2\le C\cdot\ln\left(\into n^2(\cdot,t)+1\right)+C
 	\end{align}
 	for all $t\in(0,T_{\max})$.
 \end{lem}
 \begin{proof}
 	Using mass conservation \eqref{eqn-massconservation}, we apply Lemma 2.7 of \cite{Winkler2020SIAM} to obtain $C_1>0$ such that
 	\begin{align}\label{eqn:2104-1}
 		\frac{d}{dt}\into|u|^2+\into|\nabla u|^2\le C_1\into n\ln\frac{n}{\overline{n}_0}+C_1\quad\text{for all }t\in(0,T_{\max}),
 	\end{align}
 	where $\overline{n}_0:=\frac{m}{|\Omega|}$. {Moreover, Lemma~\ref{lem:logarithmic-entropy}, applied to $\varphi=n$, yields}
 	\begin{align*}
 		\into n\ln n&\le\left(\into n+1\right)\cdot\ln\left(\into n^2+e\right)\\
 		&\le\left(\into n+1\right)\cdot\left(\ln\left(\into n^2+1\right)+1\right)\quad\text{for all }t\in(0,T_{\max}),
 	\end{align*}
 	which, together with mass conservation \eqref{eqn-massconservation}, shows that there exists $C_2>0$ such that
 	\begin{align*}
 		C_1\into n\ln\frac{n}{\overline{n}_0}=C_1\into n\ln n-C_1\ln\overline{n}_0\into n\le C_2\ln\left(\into n^2+1\right)+C_2\quad\text{for all }t\in(0,T_{\max}).
 	\end{align*}
 	Together with \eqref{eqn:2104-1}, this implies \eqref{eqn:2104-2}.
 \end{proof}
 
We next show that $\mathscr{E}(t)$ can be bounded from below provided $m<4\pi$. 
 \begin{lem}\label{lem:lowerboundofentropy}
 	Assume that \eqref{eq:newthreshold} holds. There exists $\mathcal{M}\ge 1$ such that $\mathcal{M}+\mathscr{E}(t)\ge1$ for all $t\in(0,T_{\max})$ and
 	\begin{align}\label{eqn:lowboundofentropy}
 		\mathcal{M}+\mathscr{E}(t)\ge C\cdot\left(\into n(\cdot,t)\ln\frac{n(\cdot,t)}{\overline{n}_0}+\into|\nabla c(\cdot,t)|^2\right)+\frac{1}{2}\into c^2(\cdot,t)
 	\end{align}
 	for all $t\in(0,T_{\max})$, where $C:=\frac{1}{2}-\frac{1}{2}\left(\frac{m}{4\pi}\right)^{\frac{1}{3}}>0$.
 \end{lem}
 \begin{proof}
 	{Applying Lemma~\ref{lem:entropy-interaction} with $\phi=n$, $\psi=c$, $a=\alpha:=(\frac{4\pi}{m})^{\frac{1}{3}}$ and $\varepsilon=\alpha-1>0$, and using \eqref{eqn-massconservation} and \eqref{eqn-signalmass}, we obtain $C_1>0$ such that}
 	\begin{align}\label{eqn-1404-1}
 		\into nc\le\left(\frac{m}{4\pi}\right)^{\frac{1}{3}}\into n\ln\frac{n}{\overline{n}_0}+\frac{1}{2}\left(\frac{m}{4\pi}\right)^{\frac{1}{3}}\into|\nabla c|^2+C_1\quad\text{for all }t\in(0,T_{\max}).
 	\end{align}
 	Note that $m<4\pi$ guarantees that $\left(\frac{m}{4\pi}\right)^{\frac{1}{3}}<1$, and that $\into n\ln n=\into n\ln\frac{n}{\overline{n}_0}+m\ln\frac{m}{|\Omega|}$ for all $t\in(0,T_{\max})$. We hence know from the definition of $\mathscr{E}$ in \eqref{eqn:def of entropy} and estimate \eqref{eqn-1404-1} that
 	\begin{align}\label{eqn:2234-1}
 		\mathscr{E}\ge\left(1-\left(\frac{m}{4\pi}\right)^{\frac{1}{3}}\right)\into n\ln\frac{n}{\overline{n}_0}+m\ln\frac{m}{|\Omega|}+\left(\frac{1}{2}-\frac{1}{2}\left(\frac{m}{4\pi}\right)^{\frac{1}{3}}\right)\into|\nabla c|^2+\frac{1}{2}\into c^2-C_1
 	\end{align}
 	for all $t\in(0,T_{\max})$. We now choose $\mathcal{M}\ge 1$ sufficiently large to ensure $\mathcal{M}+m\ln\frac{m}{|\Omega|}-C_1\ge 1$. Therefore, adding this $\mathcal{M}$ to both sides of \eqref{eqn:2234-1} yields \eqref{eqn:lowboundofentropy} and $\mathcal M+\mathscr E\ge1$ for all $t\in(0,T_{\max})$.
 \end{proof}
 

 \subsection{Dissipation-based control of density and signal concentration}
 
This subsection aims to establish a bound for $\int_\Omega n^2$ in terms of the two dissipation terms in \eqref{eqn:entropy}. We first apply the interpolation inequality in Lemma~\ref{lem:GNS} to $(n+k)^{1/4}-k^{1/4}$ to establish a relationship between $\into|\nabla(n+k)^{\frac{1}{4}}|^2$ and the dissipation terms.
\begin{lem}\label{lem:staticabsorbstart}
Let \eqref{eq:newthreshold} hold. There exists $C>0$ such that for every $k>0$,
 	\begin{align}\label{eq:staticabsorbstart}
 		&\int_\Om\left|\nabla(n(\cdot,t)+k)^{\frac{1}{4}}\right|^2\nonumber\\
 		&\quad\le C\cdot k^{-\frac{1}{2}}\left(\into n(\cdot,t)|\nabla(\ln n(\cdot,t)-c(\cdot,t))|^2\right)\nonumber\\
 		&\quad\quad+C\cdot\left(\into|\Delta c(\cdot,t)-c(\cdot,t)+n(\cdot,t)|^2\right)^{\frac{1}{2}}+C\cdot\sqrt{k}+C
 	\end{align}
 	for all $t\in(0,T_{\max})$.
 \end{lem}
 \begin{proof}
 	We first see that
 	\begin{align}\label{eq0946-1}
 		4\into\left|\nabla(n+k)^{\frac{1}{4}}\right|^2&=\frac{1}{4}\into\frac{|\nabla n|^2}{(n+k)^{\frac{3}{2}}}\nonumber\\
 		&=\frac{1}{4}\into\frac{\nabla n\cdot(\nabla n-n\nabla c)}{(n+k)^{\frac{3}{2}}}+\into\frac{n}{4(n+k)^{\frac{3}{2}}}\nabla n\cdot\nabla c
 	\end{align}
 	for $t\in(0,T_{\max})$. Denoting $h(s):=\frac{(\sqrt{s+k}-\sqrt k)^2}{2\sqrt{s+k}}$, $s\ge 0$, we have that
 	\begin{align*}
 		h'(s)=\frac{s}{4(s+k)^{\frac{3}{2}}}\quad\text{for }s\ge 0,
 	\end{align*}
 	and
 	\begin{equation}\label{eq:bkproperties}
 		0\le h(s)\le\frac{\sqrt s}{2}\quad\text{for }s\ge 0.
 	\end{equation}
 	Hence, integrating by parts, we deduce from \eqref{eq0946-1} that
 	\begin{align}\label{eqn:2236-1}
 		4\into\left|\nabla(n+k)^{\frac{1}{4}}\right|^2&=\frac{1}{4}\into\frac{\nabla n\cdot(\nabla n-n\nabla c)}{(n+k)^{\frac{3}{2}}}+\into\nabla h(n)\cdot\nabla c\nonumber\\
 		&=\frac{1}{4}\into\frac{\nabla n\cdot(\nabla n-n\nabla c)}{(n+k)^{\frac{3}{2}}}-\into h(n)\Delta c\nonumber\\
 		&=:I+II
 	\end{align}
 	for all $t\in(0,T_{\max})$.
 	
 	For $I$, we utilize H\"older's inequality to see that
 	\begin{align}\label{eq:Qupper}
 		\left|I\right|&=\left|\frac{1}{4}\into\frac{\nabla n\cdot(n\nabla(\ln n-c))}{(n+k)^{\frac{3}{2}}}\right|\nonumber\\
 		&\le\frac14\left(\int_\Om n|\nabla(\ln n-c)|^2\right)^{\frac{1}{2}}\cdot\left(\int_\Om\frac{n|\nabla n|^2}{(n+k)^3}\right)^{\frac{1}{2}}\notag\\
 		&\le k^{-\frac{1}{4}}\left(\int_\Om n|\nabla(\ln n-c)|^2\right)^{\frac{1}{2}}\cdot\left(\int_\Om\left|\nabla(n+k)^{\frac{1}{4}}\right|^2\right)^{\frac{1}{2}}
 	\end{align}
 	for all $t\in(0,T_{\max})$, {where} we used the fact that $\frac{n}{(n+k)^{\frac{3}{2}}}\le k^{-\frac{1}{2}}$ in the last step.
 	
 	As for $II$, we first compute that 
 	\begin{align}\label{eq1018-1}
 		II=-\into h(n)(\Delta c-c+n)-\into h(n)c+\into h(n)n\quad\text{for all }t\in(0,T_{\max}).
 	\end{align}
 	Since $c$ is nonnegative, one first has from \eqref{eq:bkproperties} that
 	\begin{align}\label{eq0953-1}
 		-\into h(n)c\le 0\quad\text{for all }t\in(0,T_{\max}),
 	\end{align}
 	and
 	\begin{align}\label{eqn:2236-1-1}
 		\into h(n)n\le\frac{1}{2}\into n^{\frac{3}{2}}\quad\text{for all }t\in(0,T_{\max}).
 	\end{align}
 	Moreover, H\"older's inequality and \eqref{eq:bkproperties} guarantee that
 	\begin{align}\label{eqn:2236-2}
 		-\into h(n)(\Delta c-c+n)&\le\|h(n)\|_{L^2(\Omega)}\|\Delta c-c+n\|_{L^2(\Omega)}\nonumber\\
 		&\le\frac{1}{2}\left(\into n\right)^{\frac{1}{2}}\cdot\|\Delta c-c+n\|_{L^2(\Omega)}\nonumber\\
 		&=\frac{\sqrt{m}}{2}\|\Delta c-c+n\|_{L^2(\Omega)}\quad\text{for all }t\in(0,T_{\max}).
 	\end{align}
 	Collecting \eqref{eq0953-1}, \eqref{eqn:2236-1-1}, and \eqref{eqn:2236-2}, we infer from \eqref{eq1018-1} that there is $C_1>0$ such that
 	\begin{align}\label{eq:1020-1}
 		II\le C_1\left(\into|\Delta c-c+n|^2\right)^{\frac{1}{2}}+\frac{1}{2}\into n^{\frac{3}{2}}\quad\text{for all }t\in(0,T_{\max}).
 	\end{align}
 	Consequently, we deduce from \eqref{eqn:2236-1}, \eqref{eq:Qupper}, and \eqref{eq:1020-1}  that there exists $C_2>0$ such that
 	\begin{align}\label{eq:1000}
 		&4\into\left|\nabla(n+k)^{\frac{1}{4}}\right|^2\nonumber\\
 		&\le k^{-\frac{1}{4}}\left(\int_\Om n|\nabla(\ln n-c)|^2\right)^{\frac{1}{2}}\cdot\left(\int_\Om\left|\nabla(n+k)^{\frac{1}{4}}\right|^2\right)^{\frac{1}{2}}+C_2\left(\into|\Delta c-c+n|^2\right)^{\frac{1}{2}}+\frac{1}{2}\into n^{\frac{3}{2}}
 	\end{align}
 	for all $t\in(0,T_{\max})$.
 	
 	We now deal with the remaining term $\frac{1}{2}\int_\Omega n^{\frac{3}{2}}$ by Lemma~\ref{lem:GNS}. Since $0<m<4\pi$, we choose
 	\[
 	\rho:=\left(\frac{4\pi}{m}\right)^{\frac14}>1,
 	\qquad
 	\varepsilon:=\frac{2}{\pi}(\rho-1)>0.
 	\]
 	Then
 	\begin{equation}\label{eq:delta}
 		4-\frac12\rho\left(\frac2\pi+\varepsilon\right)m
 		=4\left(1-\sqrt{\frac{m}{4\pi}}\right)>0.
 	\end{equation}
 	Note that for every $k>0$ and $s\ge (1-\rho^{-\frac{1}{6}})^{-4}k$, we have $s^{\frac{3}{2}}\le\rho((s+k)^{\frac14}-k^{\frac14})^6$, whereas for $0\le s<(1-\rho^{-\frac{1}{6}})^{-4}k$, one has $s^{\frac{3}{2}}\le(1-\rho^{-\frac{1}{6}})^{-2}\sqrt{k}s$. Therefore, there exists $C_3>0$, independent of $k$, such that
 	\begin{equation}\label{eq:nthreehalf}
 		\int_\Omega n^{\frac32}\le\rho\int_\Omega\left[(n+k)^{\frac14}-k^{\frac14}\right]^6+C_3\sqrt{k}\quad\text{for all }t\in(0,T_{\max}).
 	\end{equation}
 Notice that the pointwise inequality $0\le (n+k)^{\frac{1}{4}}-k^{\frac{1}{4}}\le n^{\frac{1}{4}}$ implies
 	\begin{equation*}
 		\int_\Om\left[(n+k)^{\frac{1}{4}}-k^{\frac{1}{4}}\right]^4\le \into n=m\quad\text{for all }t\in(0,T_{\max}).
 	\end{equation*}
 Applying Lemma \ref{lem:GNS} to $(n+k)^{\frac{1}{4}}-k^{\frac{1}{4}}$, we deduce from \eqref{eq:nthreehalf} that there exists {$C_4>0$, independent of $k$,} such that
 	\begin{align}\label{eq:2322-1}
 		\into n^{\frac{3}{2}}\le\rho\left(\frac2\pi+\varepsilon\right)m\int_\Om\left|\nabla(n+k)^{\frac{1}{4}}\right|^2+C_4+C_4\sqrt{k}\quad\text{for all }t\in(0,T_{\max}).
 	\end{align}
 	Therefore, combining \eqref{eq:1000} and \eqref{eq:2322-1}, we use Young's inequality to find $C_5>0$ such that
 	\begin{align}\label{eq:Qlower}
 		&4\into\left|\nabla(n+k)^{\frac{1}{4}}\right|^2\nonumber\\
 		&\le k^{-\frac{1}{4}}\left(\int_\Om n|\nabla(\ln n-c)|^2\right)^{\frac{1}{2}}\cdot\left(\int_\Om\left|\nabla(n+k)^{\frac{1}{4}}\right|^2\right)^{\frac{1}{2}}+C_2\left(\into|\Delta c-c+n|^2\right)^{\frac{1}{2}}\nonumber\\
 		&\quad+\frac{1}{2}\rho\left(\frac2\pi+\varepsilon\right)m\int_\Om\left|\nabla(n+k)^{\frac{1}{4}}\right|^2+C_4+C_4\sqrt{k}\nonumber\\
 		&\le \left(1-\sqrt{\frac{m}{4\pi}}\right)\into\left|\nabla(n+k)^{\frac{1}{4}}\right|^2+C_5k^{-\frac{1}{2}}\left(\into n|\nabla(\ln n-c)|^2\right)+C_2\left(\into|\Delta c-c+n|^2\right)^{\frac{1}{2}}\nonumber\\
 		&\quad+\frac{1}{2}\rho\left(\frac2\pi+\varepsilon\right)m\int_\Om\left|\nabla(n+k)^{\frac{1}{4}}\right|^2+C_4+C_4\sqrt{k}
 	\end{align}
 	for all $t\in(0,T_{\max})$. Combining \eqref{eq:delta} and \eqref{eq:Qlower}, we arrive at
 	\begin{align*}
 		&3\left(1-\sqrt{\frac{m}{4\pi}}\right)\into\left|\nabla(n+k)^{\frac{1}{4}}\right|^2\nonumber\\
 		&\le C_5k^{-\frac{1}{2}}\left(\into n|\nabla(\ln n-c)|^2\right)+C_2\left(\into|\Delta c-c+n|^2\right)^{\frac{1}{2}}+C_4+C_4\sqrt{k}
 	\end{align*}
 	for all $t\in(0,T_{\max})$, which yields \eqref{eq:staticabsorbstart}.
 \end{proof}
 
With the help of Lemma \ref{lem:staticabsorbstart}, we select some special $k$ to control density and signal concentration by the dissipation.
\begin{lem}\label{lem:static}
Assume that \eqref{eq:newthreshold} holds. There exists $C>0$ such that
\begin{align}\label{eq:static}
&\|n(\cdot,t)\|_{L^2(\Omega)}^2+\|c(\cdot,t)\|_{W^{2,2}(\Omega)}^2\nonumber\\
&\quad\le C\cdot\left(1+\into n(\cdot,t)|\nabla(\ln n(\cdot,t)-c(\cdot,t))|^2+\into|\Delta c(\cdot,t)-c(\cdot,t)+n(\cdot,t)|^2\right)
\end{align}
for all $t\in(0,T_{\max})$.
\end{lem}
\begin{proof}
	Let $k:=\eta\int_\Omega n^2$, where $\eta>0$ will be chosen later. We first employ Lemma \ref{lem:staticabsorbstart} to fix $C_1>0$, independent of $k$ and $\eta$, such that
	\begin{align}\label{eq:0004-1}
		\int_\Om\left|\nabla(n+k)^{\frac{1}{4}}\right|^2&\le C_1k^{-\frac{1}{2}}\left(\into n|\nabla(\ln n-c)|^2\right)+C_1\left(\into|\Delta c-c+n|^2\right)^{\frac{1}{2}}+C_1\sqrt{k}+C_1\nonumber\\
		&= C_1\eta^{-\frac{1}{2}}\left(\into n^2\right)^{-\frac{1}{2}}\cdot\left(\into n|\nabla(\ln n-c)|^2\right)+C_1\left(\into|\Delta c-c+n|^2\right)^{\frac{1}{2}}\nonumber\\
		&\quad+C_1\eta^{\frac{1}{2}}\left(\into n^2\right)^{\frac{1}{2}}+C_1
	\end{align}
	for all $t\in(0,T_{\max})$. Notice that the Gagliardo-Nirenberg inequality in our two-dimensional setting provides $C_2>0$ such that
	\begin{align}\label{eqn1030-1}
		\|f\|_{L^8(\Omega)}^8\le C_2\|\nabla f\|_{L^2(\Omega)}^4\|f\|_{L^4(\Omega)}^4+C_2\|f\|_{L^4(\Omega)}^8\quad\text{for all }f\in W^{1,2}(\Omega).
	\end{align}
	We also note that there are constants $L_1,L_2>0$, independent of $k$, such that $s^2\le L_1\left[(s+k)^{\frac{1}{4}}-k^{\frac{1}{4}}\right]^8$ when $s\ge16k$ and $s^2\le L_2ks$ when $s<16k$. Using \eqref{eqn1030-1}, the mass conservation \eqref{eqn-massconservation}, and the fact that $0\le (n+k)^{\frac{1}{4}}-k^{\frac{1}{4}}\le n^{\frac{1}{4}}$, we obtain {$C_3>0$, independent of $k$ and $\eta$,} such that
	\begin{align}\label{eq:NfromA}
		\int_\Om n^2
		&\le L_1\int_\Om \left[(n+k)^{\frac{1}{4}}-k^{\frac{1}{4}}\right]^8+L_2km\notag\\
		&\le C_3\left(\int_\Om\left|\nabla(n+k)^{\frac{1}{4}}\right|^2\right)^2+C_3\eta\into n^2+C_3\quad\text{for all }t\in(0,T_{\max}).
	\end{align}
	{Recalling that $C_1$ and $C_3$ have been fixed independently of $\eta$, we} now choose $\eta>0$ satisfying
	\begin{equation}\label{eq:eta}
		0<\eta\le\min\left\{1,\frac1{2C_3},\frac{1}{8C_1^2C_3}\right\}.
	\end{equation}
	Then, since \eqref{eq:eta} guarantees that both $C_3\eta\le\frac{1}{2}$ and $C_1\eta^{\frac{1}{2}}(2C_3)^{\frac{1}{2}}\le\frac{1}{2}$ hold, we deduce from \eqref{eq:NfromA} that
	\begin{align}\label{eqn-2311-1}
		\into n^2\le 2C_3\left(\int_\Om\left|\nabla(n+k)^{\frac{1}{4}}\right|^2\right)^2+2C_3\quad\text{for all }t\in(0,T_{\max}).
	\end{align}
	Hence
	\begin{align}\label{eq:LfromA}
		C_1\eta^{\frac{1}{2}}\left(\into n^2\right)^{\frac{1}{2}}&\le C_1\eta^{\frac{1}{2}}(2C_3)^{\frac{1}{2}}\left(\int_\Om\left|\nabla(n+k)^{\frac{1}{4}}\right|^2+1\right)\nonumber\\
		&\le\frac12\left(\int_\Om\left|\nabla(n+k)^{\frac{1}{4}}\right|^2+1\right)\quad\text{for all }t\in(0,T_{\max}).
	\end{align}
	Combining this estimate with \eqref{eq:0004-1}, we obtain $C_4>0$ such that
	\begin{align}\label{eq:Afinalstatic}
		\into\left|\nabla(n+k)^{\frac{1}{4}}\right|^2\le C_4\left(1+\left(\into|\Delta c-c+n|^2\right)^{\frac{1}{2}}+\left(\int_\Om n|\nabla(\ln n-c)|^2\right)\left(\int_\Om n^2\right)^{-\frac{1}{2}}\right)
	\end{align}
	for all $t\in(0,T_{\max})$. We now infer from \eqref{eqn-2311-1} and \eqref{eq:Afinalstatic} that there is $C_5>0$ such that
	\begin{align*}
		\left(\int_\Om n^2\right)^{\frac{1}{2}}&\le\sqrt{2C_3}\left(\int_\Om\left|\nabla(n+k)^{\frac{1}{4}}\right|^2+1\right)\nonumber\\
		&\le C_5\left(1+\left(\into|\Delta c-c+n|^2\right)^{\frac{1}{2}}+\left(\int_\Om n|\nabla(\ln n-c)|^2\right)\left(\int_\Om n^2\right)^{-\frac{1}{2}}\right)
	\end{align*}
	for all $t\in(0,T_{\max})$. Combining this with Young's inequality, we obtain that
	\begin{align*}
		\int_\Om n^2&\le C_5\left(1+\left(\into|\Delta c-c+n|^2\right)^{\frac{1}{2}}\right)\left(\int_\Om n^2\right)^{\frac{1}{2}}+C_5\int_\Om n|\nabla(\ln n-c)|^2\\
		&\le\frac12\int_\Om n^2+C_5^2\left(1+\into|\Delta c-c+n|^2\right)+C_5\int_\Om n|\nabla(\ln n-c)|^2\quad\text{for all }t\in(0,T_{\max}).
	\end{align*}
	This shows that there exists $C_6>0$ such that
	\begin{align}\label{eq1053-1}
		\into n^2\le C_6\left(1+\into n|\nabla(\ln n-c)|^2+\into |\Delta c-c+n|^2\right)\quad\text{for all }t\in(0,T_{\max}).
	\end{align}
	Finally, writing $h=\Delta c-c+n$, one can utilize the elliptic estimates to find $C_7,C_8>0$ such that
	\[
		\|c\|_{W^{2,2}(\Omega)}^2\le C_7\|n-h\|_{L^2(\Omega)}^2\le C_8\bigl(\|n\|_{L^2(\Omega)}^2+\|h\|_{L^2(\Omega)}^2\bigr)\quad\text{for all }t\in(0,T_{\max}),
	\]
	which, together with \eqref{eq1053-1}, implies \eqref{eq:static}.
\end{proof}

\subsection{Coupling the logarithm of the free energy to the fluid energy}

For convenience, we collect the two nonnegative dissipation terms in the notation
\begin{align}\label{eq:defofdissa}
	\mathscr{D}(t):=\into n(\cdot,t)|\nabla(\ln n(\cdot,t)-c(\cdot,t))|^2+\into|\Delta c(\cdot,t)-c(\cdot,t)+n(\cdot,t)|^2,\quad t\in(0,T_{\max}).
\end{align}
\begin{lem}\label{lem:evolutionofcoupledentroenergy}
	Assume that \eqref{eq:newthreshold} holds, and let $\mathcal{M}\ge 1$ be as in Lemma \ref{lem:lowerboundofentropy}. There exist $\ell>0$ and $C>0$ such that
	\begin{align}\label{eqn:1119-1}
		&\frac{d}{dt}\left(\ln(\mathcal{M}+\mathscr{E}(t))+(\ell+1)\into|u(\cdot,t)|^2\right)+\frac{1}{4}\cdot\frac{\mathscr{D}(t)}{\mathcal{M}+\mathscr{E}(t)}+\into|\nabla u(\cdot,t)|^2\nonumber\\
		&\le C\cdot\left(1+\ln(\mathcal{M}+\mathscr{E}(t))\right)
	\end{align}
	for all $t\in(0,T_{\max})$.
\end{lem}
\begin{proof}
	Letting $h:=\Delta c-c+n$ and applying the standard elliptic estimate, we find $C_1>0$ such that
	\begin{align*}
		\|c\|_{W^{2,2}(\Omega)}^2\le C_1\|n\|_{L^2(\Omega)}^2+C_1\|h\|_{L^2(\Omega)}^2\quad\text{for all }t\in(0,T_{\max}).
	\end{align*}
	Combining this estimate with Lemma \ref{lem:static} and the definition \eqref{eq:defofdissa}, we obtain $C_2>0$ such that
	\begin{align}\label{eqn:0944-1}
		\|n\|_{L^2(\Omega)}^2+\|c\|_{W^{2,2}(\Omega)}^2\le C_2\left(1+\mathscr{D}\right)\quad\text{for all }t\in(0,T_{\max}).
	\end{align}
	Furthermore, an application of the Gagliardo-Nirenberg inequality yields $C_3>0$ such that
	\begin{align*}
		\|\nabla c\|_{L^4(\Omega)}^2&\le C_3\cdot\left(\|D^2 c\|_{L^2(\Omega)}\|\nabla c\|_{L^2(\Omega)}+\|\nabla c\|_{L^2(\Omega)}^2\right)\quad\text{for all }t\in(0,T_{\max}),
	\end{align*}
	which, together with \eqref{eqn:0944-1}, provides $C_4>0$ such that
	\begin{align*}
		\|\nabla c\|_{L^4(\Omega)}^2\le C_4\|\nabla c\|_{L^2(\Omega)}\cdot\left(1+\mathscr{D}\right)^{\frac{1}{2}}+C_4\|\nabla c\|_{L^2(\Omega)}^2\quad\text{for all }t\in(0,T_{\max}).
	\end{align*}
	Therefore, by H\"older's inequality and Young's inequality, there exists $C_5>0$ such that
	\begin{align}\label{eq:1035-1}
		&\left|\int_\Om \nabla c\cdot(\nabla u\cdot\nabla c)\right|\nonumber\\
		&\le\|\nabla u\|_{L^2(\Omega)}\|\nabla c\|_{L^4(\Omega)}^2\nonumber\\
		&\le C_4\|\nabla u\|_{L^2(\Omega)}\|\nabla c\|_{L^2(\Omega)}\left(1+\mathscr{D}\right)^{\frac{1}{2}}+C_4\|\nabla u\|_{L^2(\Omega)}\|\nabla c\|_{L^2(\Omega)}^2\nonumber\\
		&\le\frac{1}{2}\mathscr{D}+C_5\cdot\left(1+\into|\nabla u|^2\right)\cdot\left(1+\into|\nabla c|^2\right)\quad\text{for all }t\in(0,T_{\max}).
	\end{align}
	Moreover, Lemma \ref{lem:lowerboundofentropy} provides $C_6>0$ and $\mathcal{M}\ge1$ such that
	\begin{align*}
		1+\into|\nabla c|^2\le C_6\cdot\left(\mathcal{M}+\mathscr{E}\right)\quad\text{for all }t\in(0,T_{\max}).
	\end{align*}
	Inserting this bound into \eqref{eq:1035-1} and using Lemma \ref{lem:evolutionofe}, we obtain $C_7>0$ such that
	\begin{align}\label{eq:1225-1}
		\frac{d}{dt}\left(\mathcal{M}+\mathscr{E}\right)+\frac{1}{2}\mathscr{D}\le C_7\cdot\left(1+\into|\nabla u|^2\right)\cdot\left(\mathcal{M}+\mathscr{E}\right)\quad\text{for all }t\in(0,T_{\max}).
	\end{align}
	We now employ Lemmas \ref{lem:2104-2} and \ref{lem:static} to find $C_8>0$ such that
	\begin{align}\label{eq:1225-2}
		(C_7+1)\frac{d}{dt}\into|u|^2+(C_7+1)\into|\nabla u|^2\le C_8\ln\left(1+\mathscr{D}\right)+C_8\quad\text{for all }t\in(0,T_{\max}).
	\end{align}
	Dividing \eqref{eq:1225-1} by $\mathcal M+\mathscr E$ and adding \eqref{eq:1225-2}, one can find $C_9>0$ such that
	\begin{align}\label{eqn:1119-2}
		\frac{d}{dt}\ln\left(\mathcal{M}+\mathscr{E}\right)+(C_7+1)\frac{d}{dt}\into|u|^2+\frac{1}{2}\frac{\mathscr{D}}{\mathcal{M}+\mathscr{E}}+\into|\nabla u|^2\le C_9+C_9\cdot\ln\left(1+\mathscr{D}\right)
	\end{align}
	for all $t\in(0,T_{\max})$. Notice that our assumption $m<4\pi$ and Lemma \ref{lem:lowerboundofentropy} ensure that $\mathcal{M}\ge 1$ and $\mathcal{M}+\mathscr{E}\ge 1$ for all $t\in(0,T_{\max})$. Due to the elementary inequality $\ln(1+\xi)\le\frac{1}{4C_9}\xi+C_{10}$ for all $\xi\ge 0$, with some $C_{10}>0$, there exists $C_{11}>0$ such that
	\begin{align*}
		C_9\ln\{1+\mathscr{D}\}&\le C_9\ln\left((\mathcal{M}+\mathscr{E})\cdot\left(1+\frac{\mathscr{D}}{\mathcal{M}+\mathscr{E}}\right)\right)\nonumber\\
		&\le C_9\ln(\mathcal{M}+\mathscr{E})+C_9\ln\left(1+\frac{\mathscr{D}}{\mathcal{M}+\mathscr{E}}\right)\nonumber\\
		&\le C_9\ln(\mathcal{M}+\mathscr{E})+\frac{1}{4}\frac{\mathscr{D}}{\mathcal{M}+\mathscr{E}}+C_{11}\quad\text{for all }t\in(0,T_{\max}).
	\end{align*}
	Combining this with \eqref{eqn:1119-2}, we see that \eqref{eqn:1119-1} is valid by letting $\ell:=C_7$.
\end{proof}
We can obtain the $t$-dependent bounds for $n\ln n$ and $\nabla c$ directly by integrating \eqref{eqn:1119-1}. To obtain bounds independent of time, we next relate the free energy to a sublinear power of the dissipation.
\begin{lem}\label{lem:upperboundofentropy}
	Assume that \eqref{eq:newthreshold} holds. Let $\mathscr{E}(t)$ and $\mathscr{D}(t)$ be defined by \eqref{eqn:def of entropy} and \eqref{eq:defofdissa}, respectively, and let $\mathcal{M}$ be as in Lemma \ref{lem:lowerboundofentropy}. There exists $C>0$ such that
	\begin{align}\label{0017-3}
		\mathcal{M}+\mathscr{E}(t)\le C\cdot\left(1+\mathscr{D}(t)\right)^{\frac{2}{3}}\quad\text{for all }t\in(0,T_{\max}).
	\end{align}
\end{lem}
\begin{proof}
	We utilize Lemma \ref{lem:static} to find $C_1>0$ such that
	\begin{align}\label{eqn-0001-1}
		\|c\|_{W^{2,2}(\Omega)}^2\le C_1\left(1+\mathscr{D}\right)\quad\text{for all }t\in(0,T_{\max}).
	\end{align}
	The Gagliardo-Nirenberg inequality provides $C_2>0$ satisfying
	\begin{align*}
		\|c\|_{L^2(\Omega)}^2\le C_2\|\Delta c\|_{L^2(\Omega)}^{\frac{2}{3}}\|c\|_{L^1(\Omega)}^{\frac{4}{3}}+C_2\|c\|_{L^1(\Omega)}^2\quad\text{for all }t\in(0,T_{\max}),
	\end{align*}
	which, together with \eqref{eqn-signalmass} and \eqref{eqn-0001-1}, yields $C_3,C_4>0$ such that
	\begin{align}\label{0017-1}
		\into c^2&\le C_3\left(\into|\Delta c|^2\right)^{\frac{1}{3}}+C_3\le C_4\left(1+\mathscr{D}\right)^{{\frac{1}{3}}}\quad\text{for all }t\in(0,T_{\max}).
	\end{align}
	Furthermore, in view of \eqref{eqn-0001-1} and \eqref{0017-1}, integrating by parts we infer from H\"older's inequality that there exists $C_5>0$ such that
	\begin{align}\label{0017-4}
		\into|\nabla c|^2&=-\into c\Delta c\nonumber\\
		&\le\|c\|_{L^2(\Omega)}\|\Delta c\|_{L^2(\Omega)}\nonumber\\
		&\le C_5\left(1+\mathscr{D}\right)^{\frac{1}{6}+\frac{1}{2}}= C_5\left(1+\mathscr{D}\right)^{\frac{2}{3}}\quad\text{for all }t\in(0,T_{\max}).
	\end{align}
	{For the density entropy, thanks to Lemma~\ref{lem:logarithmic-entropy} and mass conservation \eqref{eqn-massconservation}, we know}
	\begin{align*}
		\into n\ln n&\le\left(\into n+1\right)\cdot\ln\left(\into n^2+e\right)\le(m+1)\left(\ln\left(\into n^2+1\right)+1\right)
	\end{align*}
	for all $t\in(0,T_{\max})$. In conjunction with Lemma \ref{lem:static}, this provides $C_6>0$ fulfilling
	\begin{align*}
		\into n\ln n\le C_6\left(\ln\left(1+\mathscr{D}\right)+1\right)\quad\text{for all }t\in(0,T_{\max}).
	\end{align*}
	Since $\ln\xi\le C_7\xi^{\frac{2}{3}}$ for $\xi\ge 1$ with some $C_7>0$, there exists $C_{8}>0$ such that
	\begin{align}\label{0017-2}
		\into n\ln n\le C_8\left(1+\mathscr{D}\right)^{\frac{2}{3}}\quad\text{for all }t\in(0,T_{\max}).
	\end{align}
	Because $n,c$ are nonnegative in $\Omega\times(0,T_{\max})$, it follows from \eqref{0017-1}, \eqref{0017-4}, and \eqref{0017-2} that one can find $C_{9}>0$ such that
	\begin{align*}
		\mathcal{M}+\mathscr{E}&=\mathcal{M}+\into n\ln n+\frac{1}{2}\into|\nabla c|^2+\frac{1}{2}\into c^2-\into nc\nonumber\\
		&\le{C_{9}}\left(1+\mathscr{D}\right)^{\frac{2}{3}}
	\end{align*}
	for all $t\in(0,T_{\max})$, which leads to \eqref{0017-3}.
\end{proof}

We now combine Lemmas \ref{lem:evolutionofcoupledentroenergy} and \ref{lem:upperboundofentropy} to obtain the basic bound for $n$ and $\nabla c$, as well as  the space-time estimate for dissipation, which provide the inputs for the regularity estimates in the next subsection.
\begin{lem}\label{lem:boundofllogl}
	Assume that \eqref{eq:newthreshold} holds. There exists $C>0$ such that
	\begin{align}\label{eq0125-2}
		\into n(\cdot,t)\ln n(\cdot,t)+\into|\nabla c(\cdot,t)|^2\le C\quad\text{for all }t\in(0,T_{\max}),
	\end{align}
	and
	\begin{align}\label{eq0125-1}
		\int_t^{t+\tau}\into|\nabla u|^2\le C\quad\text{for all }t\in(0,T_{\max}-\tau),
	\end{align}
	as well as
	\begin{align}\label{eq0125-3}
		\int_t^{t+\tau}\into n|\nabla(\ln n-c)|^2+\int_t^{t+\tau}\into|\Delta c-c+n|^2\le C\quad\text{for all }t\in(0,T_{\max}-\tau),
	\end{align}
	where $\tau:=\min\{1,\frac{1}{2}T_{\max}\}$.
\end{lem}
\begin{proof}
	We first employ Lemma \ref{lem:evolutionofcoupledentroenergy} to fix some $C_1>0$ such that
	\begin{align}\label{eq1217-1}
		\frac{d}{dt}\left(\ln(\mathcal{M}+\mathscr{E})+(\ell+1)\into|u|^2\right)+\frac{1}{4}\cdot\frac{\mathscr{D}}{\mathcal{M}+\mathscr{E}}+\into|\nabla u|^2\le C_1\cdot\left(1+\ln(\mathcal{M}+\mathscr{E})\right)
	\end{align}
	for all $t\in(0,T_{\max})$. Since Lemma \ref{lem:lowerboundofentropy} provides $\mathcal{M}+\mathscr{E}\ge 1$ for all $t\in(0,T_{\max})$, we use Lemma \ref{lem:upperboundofentropy} to find $C_2>0$ satisfying
	\begin{align*}
		\left(\mathcal{M}+\mathscr{E}\right)^{\frac{1}{2}}\le C_2\cdot\frac{\mathscr{D}}{\mathcal{M}+\mathscr{E}}+C_2\quad\text{for all }t\in(0,T_{\max}).
	\end{align*}
	Note that for any $\varepsilon>0$, there exists $K>0$ such that $C_1\ln\xi\le\varepsilon\xi^{\frac{1}{2}}+K$ for all $\xi\ge 1$. It follows that one can find $C_3,C_4>0$ such that
	\begin{align}\label{0049-1}
		C_1\ln(\mathcal{M}+\mathscr{E})\le\frac{1}{16C_2}(\mathcal{M}+\mathscr{E})^{\frac{1}{2}}+C_3\le \frac{1}{16}\cdot\frac{\mathscr{D}}{\mathcal{M}+\mathscr{E}}+C_4\quad\text{for all }t\in(0,T_{\max}).
	\end{align}
	Thanks to Poincar\'e's inequality, there exists $C_5>0$ such that $(\ell+1)\into|u|^2\le C_5\into|\nabla u|^2$ for all $t\in(0,T_{\max})$. We then choose a $\delta\in(0,1)$ sufficiently small with the property that $\delta\cdot C_5\le\frac{1}{2}$ to see that
	\begin{align}\label{1416-1}
		\delta(\ell+1)\into|u|^2\le C_5\delta\into|\nabla u|^2\le\frac{1}{2}\into|\nabla u|^2\quad\text{for all }t\in(0,T_{\max}).
	\end{align}
	By \eqref{eq1217-1} and \eqref{0049-1}, we have that
	\begin{align*}
		\frac{d}{dt}\left(\ln(\mathcal{M}+\mathscr{E})+(\ell+1)\into|u|^2\right)+\frac{3}{16}\cdot\frac{\mathscr{D}}{\mathcal{M}+\mathscr{E}}+\into|\nabla u|^2\le C_1+C_4
	\end{align*}
	for all $t\in(0,T_{\max})$. Choose $a:=\min\{C_1,\delta\}>0$, where $\delta$ is as in \eqref{1416-1}. Since $\ln (\mathcal{M}+\mathscr{E})\ge0$, estimates \eqref{0049-1} and \eqref{1416-1} imply
	\begin{align*}
		a\left(\ln(\mathcal{M}+\mathscr{E})+(\ell+1)\into|u|^2\right)&\le C_1\ln(\mathcal{M}+\mathscr{E})+\delta(\ell+1)\into|u|^2\nonumber\\
		&\le\frac{1}{16}\frac{\mathscr{D}}{\mathcal{M}+\mathscr{E}}+C_4+\frac{1}{2}\into|\nabla u|^2\quad\text{for all }t\in(0,T_{\max}).
	\end{align*}
	Thus, letting $C_5:=C_1+2C_4$, $y(t):=\ln(\mathcal{M}+\mathscr{E}(t))+(\ell+1)\into|u(\cdot,t)|^2$, $t\in(0,T_{\max})$, we obtain
	\begin{align}\label{0058-1}
		y'(t)+ay(t)+\frac18\cdot\frac{\mathscr{D}}{\mathcal{M}+\mathscr{E}}+\frac12\into|\nabla u|^2\le C_5\quad\text{for all }t\in(0,T_{\max}).
	\end{align}
Using the comparison principle, we can find $C_6>0$ such that
\begin{align*}
	y(t)\le e^{-at}y(0)+\frac{C_5}{a}(1-e^{-at})\le C_6\quad\text{for all }t\in(0,T_{\max}),
\end{align*}
which implies
\begin{align}\label{eq-1226-1}
\ln(\mathcal{M}+\mathscr{E})+\int_\Omega|u|^2\le C_6\quad\text{for all }t\in(0,T_{\max}).
\end{align}
Integrating \eqref{0058-1} on each interval $[t,t+\tau]$, where $\tau:=\min\{1,\frac{1}{2}T_{\max}\}$, and using \eqref{eq-1226-1}, $y\ge0$, $a>0$ and $\tau\le1$, we find $C_7>0$ such that
\begin{align}\label{0500-1}
\int_t^{t+\tau}\frac{\mathscr D}{\mathcal M+\mathscr E}
+\int_t^{t+\tau}\int_\Omega|\nabla u|^2
\le C_7\quad\text{for all }t\in(0,T_{\max}-\tau).
\end{align}
In particular, \eqref{0500-1} implies \eqref{eq0125-1}. Moreover,
\begin{align}\label{eqn1242}
1\le \mathcal{M}+\mathscr{E}= e^{\ln(\mathcal{M}+\mathscr{E})}\le C_8:=e^{C_6}\quad\text{for all }t\in(0,T_{\max}),
\end{align}
and hence
\[
\int_t^{t+\tau}\mathscr D(s)\,\mathrm ds=\int_t^{t+\tau}\frac{\mathscr{D}}{\mathcal{M}+\mathscr{E}}\cdot(\mathcal{M}+\mathscr{E})\,\mathrm ds\le C_8C_7\quad\text{for all }t\in(0,T_{\max}-\tau),
\]
which proves \eqref{eq0125-3}. Finally, since Lemma~\ref{lem:lowerboundofentropy} provides $C_8>0$ such that
\[
\int_\Omega n\ln\frac{n}{\overline n_0}+\int_\Omega|\nabla c|^2\le C_8(\mathcal{M}+\mathscr{E})\quad\text{for all }t\in(0,T_{\max}),
\]
we can deduce from \eqref{eqn1242} and the mass conservation \eqref{eqn-massconservation} that there exists $C_9>0$ such that
\[
\int_\Omega n\ln n+\int_\Omega|\nabla c|^2\le C_9\quad\text{for all }t\in(0,T_{\max}),
\]
which yields \eqref{eq0125-2}.
\end{proof}


\subsection{Higher regularity and proof of Theorem \ref{thm:main}}\label{sec:continuation}

We shall make full use of the basic regularity information obtained in Lemma \ref{lem:boundofllogl} to derive higher regularity in this subsection. We first follow the argument of \cite[Lemma~3.7]{WangWinkler2024} to estimate $\|n\|_{L^2(\Omega)}$ and $\|\nabla c\|_{L^4(\Omega)}$. For the sake of completeness, we reproduce the proof here.

\begin{lem}\label{lem:l2estimatesforn}
	Assume that \eqref{eq:newthreshold} holds. There exists $C>0$ such that
	\begin{align}\label{0207-1}
		\into n^2(\cdot,t)+\into|\nabla c(\cdot,t)|^4\le C\quad\text{for all }t\in(0,T_{\max}).
	\end{align}
\end{lem}
\begin{proof}

Note the identity that $\nabla c\cdot\nabla\Delta c=\frac12\Delta|\nabla c|^2-|D^2c|^2$ with $D^2c$ being the Hessian matrix of $c$. Integrating by parts in the first two equations of \eqref{eq:KSNS}, and using $\nabla\!\cdot u=0$, we obtain
\begin{align*}
 &\frac{d}{dt}\left(\into n^2+\into|\nabla c|^4\right)\\
 &\quad=-2\into|\nabla n|^2+2\into n\nabla n\cdot\nabla c
       +4\into|\nabla c|^2\nabla c\cdot\nabla(\Delta c-c+n-u\cdot\nabla c)\\
 &\quad=-2\into|\nabla n|^2+2\into n\nabla n\cdot\nabla c
       +2\into|\nabla c|^2\Delta|\nabla c|^2\\
 &\qquad\quad-4\into|\nabla c|^2|D^2c|^2-4\into|\nabla c|^4
       +4\into|\nabla c|^2\nabla n\cdot\nabla c\\
 &\qquad\quad-4\into|\nabla c|^2\nabla c\cdot(\nabla u\cdot\nabla c)\\
 &\quad=-2\into|\nabla n|^2+2\into n\nabla n\cdot\nabla c
       -2\into\bigl|\nabla|\nabla c|^2\bigr|^2\\
 &\qquad\quad+2\int_{\partial\Omega}|\nabla c|^2\frac{\partial|\nabla c|^2}{\partial\nu}
       -4\into|\nabla c|^2|D^2c|^2-4\into|\nabla c|^4\\
 &\qquad\quad-4\into n\left(2\nabla c\cdot(D^2c\cdot\nabla c)+|\nabla c|^2\Delta c\right)
       -4\into|\nabla c|^2\nabla c\cdot(\nabla u\cdot\nabla c)
\end{align*}
for all $t\in(0,T_{\max})$. Here the remaining transport term vanishes because
\[
 \into|\nabla c|^2\nabla c\cdot(D^2c\cdot u)
 =\frac14\into u\cdot\nabla|\nabla c|^4=0.
\]
Since $|\Delta c|\le\sqrt2\,|D^2c|$, using Young's inequality, with $C_2:=1+2(2+\sqrt2)^2$, we have
\begin{align}\label{eq1114-1}
 &\frac{d}{dt}\left(\into n^2+\into|\nabla c|^4\right)\nonumber\\
 &\quad\le-\into|\nabla n|^2+\into n^2|\nabla c|^2
       -2\into\bigl|\nabla|\nabla c|^2\bigr|^2\nonumber\\
 &\qquad+2\int_{\partial\Omega}|\nabla c|^2\frac{\partial|\nabla c|^2}{\partial\nu}
       -4\into|\nabla c|^2|D^2c|^2\nonumber\\
 &\qquad+4(2+\sqrt2)\into n|\nabla c|^2|D^2c|
       +4\into|\nabla u|\,|\nabla c|^4\nonumber\\
 &\quad\le-\into|\nabla n|^2
       -2\into\bigl|\nabla|\nabla c|^2\bigr|^2-2\into|\nabla c|^2|D^2c|^2\nonumber\\
 &\qquad+2\int_{\partial\Omega}|\nabla c|^2\frac{\partial|\nabla c|^2}{\partial\nu}+C_2\into n^2|\nabla c|^2+4\into|\nabla u||\nabla c|^4\nonumber\\
       &\quad=:-\into|\nabla n|^2
       -2\into\bigl|\nabla|\nabla c|^2\bigr|^2-2\into|\nabla c|^2|D^2c|^2+I+II+III
\end{align}
for all $t\in(0,T_{\max})$.

We now estimate $I,II,III$ one by one. The smoothness of $\partial\Omega$ ensures the existence of $C_3>0$ such that
\[
 \frac{\partial|\nabla\varphi|^2}{\partial\nu}
 \le C_3|\nabla\varphi|^2\quad\text{on }\partial\Omega
\]
for every $\varphi\in C^2(\overline\Omega)$ satisfying $\partial_\nu\varphi=0$ on $\partial\Omega$. By continuity of the trace embedding $W^{\frac{3}{4},2}(\Omega)\hookrightarrow L^2(\partial\Omega)$ and compactness of $W^{1,2}(\Omega)\hookrightarrow W^{\frac34,2}(\Omega)$, the uniform bound for $\into|\nabla c|^2$ in Lemma~\ref{lem:boundofllogl} yields $C_4>0$ such that
\begin{align}\label{eq1114-2}
	I\le\frac13\into\bigl|\nabla|\nabla c|^2\bigr|^2+C_4
	\quad\text{for all }t\in(0,T_{\max}).
\end{align}
Next, the Gagliardo-Nirenberg inequality, Lemma~\ref{lem:boundofllogl} and Young's inequality provide $C_5,C_6,C_7,C_8>0$ such that
\begin{align}\label{eq1111-1}
 II
 &\le C_2\|n\|_{L^3(\Omega)}^2\bigl\||\nabla c|^2\bigr\|_{L^3(\Omega)}\nonumber\\
 &\le C_5\|n\|_{L^3(\Omega)}^2
   \left(\bigl\|\nabla|\nabla c|^2\bigr\|_{L^2(\Omega)}^2
   \bigl\||\nabla c|^2\bigr\|_{L^1(\Omega)}
   +\bigl\||\nabla c|^2\bigr\|_{L^1(\Omega)}^3\right)^{1/3}\nonumber\\
 &\le C_6\|n\|_{L^3(\Omega)}^2
   \left(\bigl\|\nabla|\nabla c|^2\bigr\|_{L^2(\Omega)}^2+1\right)^{1/3}\nonumber\\
 &\le\frac13\left(\into\bigl|\nabla|\nabla c|^2\bigr|^2+1\right)
   +C_7\into n^3\quad\text{for all }t\in(0,T_{\max}).
\end{align}
Since Lemma~\ref{lem:boundofllogl} and mass conservation \eqref{eqn-massconservation} provide $K>0$ such that
\[
\into(n+1)\ln(n+1+e)\le K
\quad\text{for all }t\in(0,T_{\max}),
\]
we may therefore apply Lemma~\ref{lem:orlicz-interpolation} with $\varphi=n+1$, $\hbar(s)=\ln(s+e)$ to find $C_8>0$ such that
\begin{align*}
	C_7\into n^3\le\frac{1}{2}\into|\nabla n|^2+C_8\quad\text{for all }t\in(0,T_{\max}),
\end{align*}
which, combined with \eqref{eq1111-1}, shows that
\begin{align}\label{eq1114-3}
	II\le\frac13\into\bigl|\nabla|\nabla c|^2\bigr|^2+\frac12\into|\nabla n|^2+C_8\quad\text{for all }t\in(0,T_{\max}).
\end{align}
An application of the Gagliardo-Nirenberg inequality and Young's inequality entails $C_9,C_{10}>0$ such that
\begin{align}\label{eq1114-4}
 III
 &\le4\|\nabla u\|_{L^2(\Omega)}\bigl\||\nabla c|^2\bigr\|_{L^4(\Omega)}^2\nonumber\\
 &\le C_9\|\nabla u\|_{L^2(\Omega)}
   \bigl\|\nabla|\nabla c|^2\bigr\|_{L^2(\Omega)}
   \bigl\||\nabla c|^2\bigr\|_{L^2(\Omega)}\nonumber\\
 &\qquad+C_9\|\nabla u\|_{L^2(\Omega)}
   \bigl\||\nabla c|^2\bigr\|_{L^2(\Omega)}^2\nonumber\\
 &\le\frac13\into\bigl|\nabla|\nabla c|^2\bigr|^2
   +C_{10}\left(\|\nabla u\|_{L^2(\Omega)}^2+1\right)\into|\nabla c|^4
 \quad\text{for all }t\in(0,T_{\max}).
\end{align}
Combining \eqref{eq1114-1}, \eqref{eq1114-2}, \eqref{eq1114-3}, and \eqref{eq1114-4}, we obtain $C_{11}>0$ satisfying that
\begin{align}\label{0200-1}
 \frac{d}{dt}\left(\into n^2+\into|\nabla c|^4\right)
 \le C_{11}+C_{11}\cdot\left(\into|\nabla u|^2+1\right)
 \cdot\left(\into n^2+\into|\nabla c|^4\right)
\end{align}
for all $t\in(0,T_{\max})$. Lemma \ref{lem:static} and estimate \eqref{eq0125-3} also provide {$C_{12}>0$} such that
\begin{align*}
 \int_t^{t+\tau}\into n^2+\int_t^{t+\tau}\|c\|_{W^{2,2}(\Omega)}^2\le{C_{12}}\quad\text{for all }t\in(0,T_{\max}-\tau)
\end{align*}
with $\tau:=\min\{1,\frac{1}{2}T_{\max}\}$. Note that the Gagliardo-Nirenberg inequality and \eqref{eq0125-2} in Lemma \ref{lem:boundofllogl} provide {$C_{13},C_{14}>0$} such that
\begin{align*}
 \|\nabla c\|_{L^4(\Omega)}^4&\le{C_{13}}\|\nabla c\|_{L^2(\Omega)}^2\|D^2c\|_{L^2(\Omega)}^2+{C_{13}}\|\nabla c\|_{L^2(\Omega)}^4\nonumber\\
 &\le{C_{14}}\|c\|_{W^{2,2}(\Omega)}^2+{C_{14}}\quad\text{for all }t\in(0,T_{\max}).
\end{align*}
Consequently, given any $t_0\in(\tau,T_{\max})$, one can choose $t^\star=t^\star(t_0)\in(t_0-\tau,t_0)$ such that there is $C_{15}>0$ fulfilling
\begin{align*}
 \into n^2(\cdot,t^\star)+\into|\nabla c(\cdot,t^\star)|^4\le{C_{15}}.
\end{align*}
Setting $y(t):=\into n^2(\cdot,t)+\into|\nabla c(\cdot,t)|^4+1$ and $h(t):=\into|\nabla u(\cdot,t)|^2+1$ for $t\in(0,T_{\max})$, we see from \eqref{0200-1} that
\[
	y'(t)\le C_{11} h(t)y(t)\quad\text{for all }t\in(0,T_{\max}).
\]
Integrating this and using \eqref{eq0125-1}, we obtain {$C_{16}>0$} such that
\begin{align*}
 y(t_0)\le y(t^\star)e^{{C_{11}}\int_{t^\star}^{t_0}h(s)ds}\le{C_{16}}.
\end{align*}
As $t_0\in(\tau,T_{\max})$ was arbitrary, in view of the continuity of $n$ and $\nabla c$ in $\Omega\times[0,\tau]$ guaranteed by Lemma \ref{lem:localexistence}, this leads to \eqref{0207-1}.
\end{proof}

{In light of Lemma \ref{lem:forced-fluid-regularity}, this particularly implies two convenient boundedness properties of $u$.}
\begin{prop}\label{lem:boundednessforu}
	Assume that \eqref{eq:newthreshold} holds. For each $\alpha\in(\frac{1}{2},1)$, there exists $K>0$ such that
	\begin{align}\label{eq2008-1}
		\|\mathcal{A}^\alpha u(\cdot,t)\|_{L^2(\Omega)}\le K\quad\text{for all }t\in(0,T_{\max}).
	\end{align}
	In particular, there is $C>0$ such that
	\begin{align}\label{0244-2}
		\|u(\cdot,t)\|_{L^\infty(\Omega)}\le C\quad\text{for all }t\in(0,T_{\max}).
	\end{align}
\end{prop}
\begin{proof}
	{Lemma~\ref{lem:l2estimatesforn} and $\nabla\phi\in L^\infty(\Omega)$ yield a uniform $L^2$ bound for $f=n\nabla\phi$. Since $u_0\in D(\mathcal{A})$, we may apply Lemma~\ref{lem:forced-fluid-regularity} with $v=u$ and $\beta=\alpha$ on every finite interval $(0,T)$ with $T<T_{\max}$. The independence of its constant from $T$ shows that, for every $\alpha\in(\frac12,1)$, there exists $C_1>0$ such that}
	\begin{align}\label{eq2009-1}
		\|\mathcal{A}^\alpha u\|_{L^2(\Omega)}\le C_1\quad\text{for all }t\in(0,T_{\max}),
	\end{align}
	which implies \eqref{eq2008-1}. According to the embedding $D(\mathcal{A}^\alpha)\hookrightarrow L^\infty(\Omega;\mathbb{R}^2)$ for $\alpha>\frac{1}{2}$, we deduce from \eqref{eq2009-1} that \eqref{0244-2} is valid.
\end{proof}

Heat semigroup estimates then improve the spatial integrability of the signal gradient.
\begin{prop}\label{lem:lqforc}
	Assume that \eqref{eq:newthreshold} holds. For every $q>2$, there is $C>0$ such that
	\begin{align}\label{eq:lqforc}
		\|\nabla c(\cdot,t)\|_{L^q(\Omega)}\le C\quad\text{for all }t\in(0,T_{\max}).
	\end{align}
\end{prop}
\begin{proof}
	Fix $q\in(2,\infty)$. We infer from the regularization estimates of the Neumann heat semigroup $(e^{t\Delta})_{t\ge0}$ in $\Omega$ that there is $C_1>0$ such that
	\begin{align*}
		&\|\nabla c(\cdot,t)\|_{L^q(\Omega)}\nonumber \\
		&= \left\|{\nabla e^{t(\Delta_N - 1)} c_0} + \int_{0}^t{\nabla e^{(t-s)(\Delta_N-1)}} \{ n(\cdot, s) - u(\cdot, s) \cdot \nabla c(\cdot, s) \} ds \right\|_{L^q(\Omega)} \nonumber\\
		&\le C_1 \|c_0\|_{W^{1,\infty}(\Omega)}\nonumber\\
		&\quad+ C_1 \int_{0}^te^{-(t-s)} \left(1+(t-s)^{-1+\frac{1}{q}}\right) \{ \|n(\cdot, s)\|_{L^2(\Omega)} + \|u(\cdot, s)\|_{L^\infty(\Omega)} \|\nabla c\|_{L^2(\Omega)} \} ds
	\end{align*}
	for all $t\in(0,T_{\max})$. Estimate \eqref{eq:lqforc} now follows from Lemmas \ref{lem:boundofllogl}, \ref{lem:l2estimatesforn}, and Proposition \ref{lem:boundednessforu}.
\end{proof}

{We next test the density equation directly to obtain a uniform $L^4$ bound, which supplies the input for the final semigroup smoothing argument.}
\begin{prop}\label{lem:boundednessforncu}
	Assume that \eqref{eq:newthreshold} holds. There exists $C>0$ such that
	\begin{align}\label{boundednessforncu}
		\|n(\cdot,t)\|_{L^\infty(\Omega)}\le C\quad\text{for all }t\in(0,T_{\max}).
	\end{align}
\end{prop}
\begin{proof}
	Testing the first equation in \eqref{eq:KSNS} by $n^3$ and applying Young's inequality and the Gagliardo-Nirenberg inequality, we take $C_i (i=1,2,3,4)$ such that
	\begin{align*}
		\frac{1}{4} \frac{d}{dt} \int_\Omega n^4 + \frac{3}{4} \int_\Omega |\nabla n^2|^2 &= 3 \int_\Omega n^3 \nabla n \cdot \nabla c \nonumber\\
		&\le \frac{3}{16} \int_\Omega |\nabla n^2|^2 + C_1\int_\Omega n^4 |\nabla c|^2 \nonumber\\
		&\le \frac{3}{16} \int_\Omega |\nabla n^2|^2 + C_1\|\nabla c\|_{L^4(\Omega)}^2 \|n^2\|_{L^4(\Omega)}^2 \nonumber\\
		&\le \frac{3}{16} \int_\Omega |\nabla n^2|^2 + C_2 \|\nabla c\|_{L^4(\Omega)}^2 \cdot \left( \|\nabla n^2\|_{L^2(\Omega)}^2 \|n^2\|_{L^{\frac{1}{2}}(\Omega)}^{\frac{2}{7}} + \|n^2\|_{L^{\frac{1}{2}}(\Omega)}^{\frac{16}{7}} \right)^{\frac{7}{8}} \nonumber\\
		&\le \frac{3}{16} \int_\Omega |\nabla n^2|^2 + C_3 \|\nabla c\|_{L^4(\Omega)}^2 \cdot \left( \|\nabla n^2\|_{L^2(\Omega)}^2 + 1 \right)^{\frac{7}{8}} \nonumber\\
		&\le \frac{3}{16} \int_\Omega |\nabla n^2|^2 + \frac{3}{16} \cdot \left( \|\nabla n^2\|_{L^2(\Omega)}^2 + 1 \right) + C_4 \|\nabla c\|_{L^4(\Omega)}^{16}
	\end{align*}
	for all $t\in(0,T_{\max})$. Combined with Proposition \ref{lem:lqforc}, this shows that there is {$C_5>0$} satisfying
	\begin{align}\label{eqn1535-1}
		\frac{d}{dt}\into n^4+\into|\nabla n^2|^2\le C_5\quad\text{for all }t\in(0,T_{\max}).
	\end{align}
	Note that the Gagliardo-Nirenberg inequality and Lemma \ref{lem:l2estimatesforn} {provide} $C_6,C_7>0$ such that
	\begin{align*}
		\|n^2\|_{L^2(\Omega)}^2\le C_6\|\nabla n^2\|_{L^2(\Omega)}^2+C\|n^2\|_{L^1(\Omega)}^2\le C_7\|\nabla n^2\|_{L^2(\Omega)}^2+C_7\quad\text{for all }t\in(0,T_{\max}).
	\end{align*}
	Selecting some $\delta\in(0,1)$ suitably small to ensure $\delta C_7\le\frac{1}{2}$, we deduce from \eqref{eqn1535-1} that one can find $C_8>0$ fulfilling
	\begin{align*}
		\frac{d}{dt}\into n^4+\delta\into n^4\le C_8\quad\text{for all }t\in(0,T_{\max}).
	\end{align*}
	Consequently, an application of {the comparison principle allows} us to obtain $C_9>0$ such that
	\begin{align}\label{0337-1}
		\into n^4\le C_9\quad\text{for all }t\in(0,T_{\max}).
	\end{align}
	
	We now take $p\in(2,4)$ and invoke standard smoothing estimates for the Neumann heat semigroup $(e^{t\Delta})_{t\ge 0}$ to find $\lambda_1>0$ and $C_{10}>0$ such that
	\begin{align*}
		&\|n(\cdot,t)\|_{L^\infty(\Omega)}\nonumber\\
		&\le C_{10}\|n_0\|_{L^\infty(\Omega)}\nonumber\\
		&\quad+C_{10}\int_0^te^{-\lambda_1(t-s)}\left(1+(t-s)^{-\frac{1}{2}-\frac{1}{p}}\right)\left(\|n(\cdot,s)u(\cdot,s)\|_{L^p(\Omega)}+\|n(\cdot,s)\nabla c(\cdot,s)\|_{L^p(\Omega)}\right)ds\nonumber\\
		&\le C_{10}\|n_0\|_{L^\infty(\Omega)}\nonumber\\
		&\quad+C_{10}\int_0^te^{-\lambda_1(t-s)}\left(1+(t-s)^{-\frac{1}{2}-\frac{1}{p}}\right)\|n(\cdot,s)\|_{L^4(\Omega)}\cdot\left(\|\nabla c(\cdot,s)\|_{L^{\frac{4p}{4-p}}(\Omega)}+\|u(\cdot,s)\|_{L^{\frac{4p}{4-p}}(\Omega)}\right)ds
	\end{align*}
	for all $t\in(0,T_{\max})$. Combining this with \eqref{0337-1} and Propositions \ref{lem:boundednessforu} and \ref{lem:lqforc} yields \eqref{boundednessforncu}.
\end{proof}

Our overall conclusion thereby becomes evident:
\begin{proof}[Proof of Theorem \ref{thm:main}]
	By virtue of Propositions \ref{lem:boundednessforu}, \ref{lem:lqforc}, and \ref{lem:boundednessforncu}, we infer from \eqref{eqn:extensibility} that $T_{\max}=\infty$. Consequently, the solution constructed in Lemma \ref{lem:localexistence} solves \eqref{eq:KSNS} classically in $\Omega\times(0,\infty)$ and satisfies the uniform bound  \eqref{eq:uniform-bounds}.
\end{proof}

\noindent
{\bf Acknowledgments.} L. Du is supported by National Natural Science Foundation of China under Grants 12526202, 12671245, 12631009. Y. Wang is supported by National Natural Science Foundation of China (Nos. 12371210, 12471195) and Natural Science Foundation of Sichuan Province (No. 26QNJJA0020).

\noindent
{\bf AI declarations.}  OpenAI’s GPT was used to assist with language editing, consistency checking, and the organization of the manuscript. The authors are responsible for the mathematical content and conclusions of the article.


\appendix
\makeatletter
\def\@seccntformat#1{%
	\ifnum\pdfstrcmp{#1}{section}=0
	Appendix~\csname the#1\endcsname:\quad
	\else
	\csname the#1\endcsname\quad
	\fi}
\makeatother
\section{Basic properties of decreasing rearrangement}\label{sec:rearrangement}

We collect the truncation and coarea properties used in the proof of Lemmas \ref{lem:propertyofrearragement} and \ref{lem:GNS}. The proof in this section is based on Rakotoson's comprehensive monograph \emph{Relative Rearrangement}~\cite{Rak25}, which provides a systematic treatment of decreasing and relative rearrangements and their applications to partial differential equations. Certain facts from geometric measure theory (\cite{EG15,Maggi2012,BFR26}) are also invoked in the proof. 

Let $\Omega\subset\mathbb{R}^2$ be a bounded domain with smooth boundary. For a nonnegative function $0\not\equiv f\in C^1(\overline{\Omega})$ and $s_0\in(0,|\Omega|)$, we let
\begin{align*}
	b:=\left(\frac{2}{s_0}\into f^4\right)^{\frac{1}{4}},\quad g:=(f-b)_+,\quad S:=|\{x\in\Omega:g>0\}|.
\end{align*}
Following \cite[Definition 1.1.3]{Rak25}, we define the decreasing rearrangement by
\begin{align*}
	g^*(s):=\inf\{r\ge 0:|\{x\in\Omega:g>r\}|\le s\},\quad s\in(0,|\Omega|).
\end{align*}
We also define the Schwarz symmetrization of $g^*$ by
\[
G(x):=
\begin{cases}
	g^*(\pi|x|^2),&0<\pi|x|^2<S,\\
	0,&\pi|x|^2\ge S\text{ or }x=0.
\end{cases}
\]

The following lemma collects properties that hold for every \(s_0\in(0,|\Omega|)\). Parts~(i)--(iii) provide bounds on the areas of the superlevel sets and level layers of \(g\), together with the coarea identity involving the relative perimeter. Part~(iv) establishes the local regularity of \(g^*\), its vanishing trace at \(s=S\) when \(S>0\), and the Sobolev membership and radial gradient identity for \(G\). In Lemma~\ref{lem:propertyofrearragement}, \(s_0\) is chosen according to the small-volume estimate in Lemma~\ref{lem:reliso}; the level-set bounds, the coarea identity, and the local regularity of \(g^*\) then yield the weighted estimate~\eqref{eq:1558-1}. Combining this estimate with part~(iv) gives the gradient comparison~\eqref{eq:radialgradient} and justifies the application of the sharp planar inequality~\eqref{eq-gns} to \(G\) in the proof of Lemma~\ref{lem:GNS}.

\begin{lem}\label{lem:some basic properties of the rearrangement}
	The following properties hold.
	\begin{itemize}
		\item[(i)] $g\in W^{1,\infty}(\Omega)$, $\into|\nabla g|^2\le\into|\nabla f|^2$, and $S\le\frac{s_0}{2}$;

		\item[(ii)] For $0<s<t\le S$, we have
		\begin{align}\label{eqn:1451-1}
			|\{x\in\Omega:g^*(t)<g<g^*(s)\}|\le t-s,
		\end{align}
		and
		\begin{align}\label{eqn:1451-2}
			s<|\{x\in\Omega:g>r\}|\le t\le s_0\quad\text{for }g^*(t)<r<g^*(s);
		\end{align}
		
		\item[(iii)] For $0\le\alpha<\beta<\infty$,
		\begin{align}\label{eqn:1522-1}
			\int_\alpha^\beta\mathscr{P}_\Omega(\{x\in\Omega:g>r\})dr=\int_{\{x\in\Omega:\alpha<g<\beta\}}|\nabla g|.
		\end{align}

		\item[{(iv)}] $G\in W^{1,2}(\mathbb R^2)$. Moreover, if $S>0$,
        $g^*$ is Lipschitz on every $[\delta,S]$ for $0<\delta<S$. We also have $\lim_{s\to S^-}g^*(s)=g^*(S)=0$, and
        \begin{align}\label{a.4}
        	\int_{\mathbb R^2}|\nabla G|^2dx=4\pi\int_0^S s\left|(g^*)'(s)\right|^2ds.
        \end{align}
	\end{itemize}
\end{lem}
\begin{proof}
		{\it Proof of (i).} By \cite[Theorem 4.4]{EG15} and the assumption $f\in C^1(\overline{\Omega})$, we have $g\in W^{1,\infty}(\Omega)$ and $\into|\nabla g|^2\le\into|\nabla f|^2$. Chebyshev's inequality yields $S=|\{x\in\Omega:g>0\}|\le\frac{1}{b^4}\into f^4=\frac{s_0}{2}$.
	
	\textit{Proof of (ii).} Let
\[
\mu(r):=|\{x\in\Omega:g(x)>r\}|,\qquad r\ge0.
\]
By \cite[Property 1.1.1(c)]{Rak25},
\begin{equation*}
\mu(g^*(q))\le q\qquad\text{for }0<q<|\Omega|.
\end{equation*}
On the other hand, we infer from the definition of $g^*$ that $\mu(r)>q$ whenever $0\le r<g^*(q)$. If $g^*(q)>0$, letting $r\nearrow g^*(q)$ and using continuity of measure from above therefore yields
\[
|\{x\in\Omega:g(x)\ge g^*(q)\}|\ge q.
\]
Consequently,
\begin{equation}\label{eq2303-1}
|\{g>g^*(q)\}|\le q\le|\{g\ge g^*(q)\}|.
\end{equation}
Now fix $0<s<t\le S$. We first deduce from part (i) that $t\le S\le s_0/2<|\Omega|$. Since $g^*$ is nonincreasing, $g^*(t)\le g^*(s)$, which implies that $\{g\ge g^*(s)\}\subset\{g>g^*(t)\}$. Hence, we infer from \eqref{eq2303-1} that
\begin{align*}
|\{g^*(t)<g<g^*(s)\}|
&=|\{g>g^*(t)\}|-|\{g\ge g^*(s)\}|\\
&\le t-s,
\end{align*}
which proves \eqref{eqn:1451-1}. Finally, if $g^*(t)<r<g^*(s)$, we have from the definition of $g^*(s)$ and the monotonicity of $\mu$ that
\[
s<\mu(r)\le\mu(g^*(t))\le t\le S\le s_0/2\le s_0.
\]
This shows \eqref{eqn:1451-2}.
	
	{\it Proof of (iii).} Let $d\equiv d(x):=\text{dist}(x,\mathbb{R}^2\setminus\Omega)$ and set $\Omega_j:=\{x\in\Omega:d>\frac{1}{j}\}$. We then choose a family of cut-off functions $(\chi_j)_{j\in\mathbb{N}}\subset C_0^\infty(\Omega)$ satisfying
	\begin{align*}
		0\le\chi_j\le 1,\quad\chi_j=1\text{ if }x\in\left\{x\in\Omega:d>\frac{3}{4j}\right\},\quad\text{supp}\chi_j\subset\left\{x\in\Omega:d\ge\frac{1}{4j}\right\}.
	\end{align*}
	For any fixed $0\le\alpha<\beta<\infty$, we define
	\begin{align*}
		\mathfrak{g}_j(x):=\left\{\begin{aligned}
			&\chi_j(x)\min\left((g(x)-\alpha)_+,\beta-\alpha\right), &&x\in\Omega,\\
			&0, &&x\in\mathbb{R}^2\setminus\Omega.
		\end{aligned}\right.
	\end{align*}
	Applying the Fleming-Rishel coarea formula (see \cite[Theorem 13.1]{Maggi2012} or \cite[Theorem 5.1.14]{BFR26}), we have
	\begin{align}\label{eqn:3.11}
		\int_{\Omega_j}|\nabla\mathfrak{g}_j|=\int_{\mathbb{R}}\mathscr{P}_{\Omega_j}(\{x\in\Omega:\mathfrak{g}_j>t\})dt.
	\end{align}
	For $0<t<\beta-\alpha$, $\{x\in\Omega:\mathfrak{g}_j>t\}\cap\Omega_j=\{x\in\Omega:g>\alpha+t\}\cap\Omega_j$ and $\nabla\mathfrak{g}_j=\textbf{1}_{\{x\in\Omega:\alpha<g<\beta\}}\nabla g$ a.e. in $\Omega_j$. Substituting these identities into \eqref{eqn:3.11}, we have
	\begin{align}\label{eq:1521-1}
		\int_{\Omega_j\cap\{x\in\Omega:\alpha<g<\beta\}}|\nabla g|=\int_\alpha^\beta\mathscr{P}_{\Omega_j}(\{x\in\Omega:g>r\})dr.
	\end{align}
	For every measurable $E$, we have $\mathscr{P}_{\Omega_j}(E)\nearrow\mathscr{P}_\Omega(E)$ as $j\to\infty$. Monotone convergence in \eqref{eq:1521-1} therefore yields \eqref{eqn:1522-1}.

\textit{Proof of (iv).}
If $S=0$, then $G=0$ almost everywhere, and the assertion
is immediate. Assume $S>0$. Since $\Omega$ is connected and has smooth boundary, the relative isoperimetric inequality \cite[Lemma~3.4.1(2), p.~82]{Rak25} provides a constant $Q>0$ such that
\begin{align}\label{eq:2033-1}
	\min\left\{|E|^{\frac{1}{2}},|\Omega\setminus E|^{\frac{1}{2}}\right\}\le Q\,\mathscr P_\Omega(E)
\end{align}
for every measurable set $E\subset\Omega$ of finite relative perimeter. Fix
\[
	\mathscr{C}:=\min\left\{\pi,Q^{-2}\right\}\in(0,2\pi),
\]
which, together with \eqref{eq:2033-1}, yields that whenever $|E|\le\frac{|\Omega|}{2}$,
\begin{align}\label{eq2313-1}
	\mathscr P_\Omega(E)\ge\sqrt{\mathscr{C}|E|}.
\end{align}
Let $0<s<t\le S$. If $g^*(s)=g^*(t)$, then \eqref{2142-1} below is immediate. Otherwise, part~(ii), specifically \eqref{eqn:1451-2}, and part~(i) show that
\begin{align}\label{eq2313-2}
	s<|\{g>r\}|\le t\le S\le\frac{s_0}{2}<\frac{|\Omega|}{2}\quad
	\text{for }g^*(t)<r<g^*(s).
\end{align}
The coarea identity in part~(iii) ensures that these superlevel sets have finite relative perimeter for almost every such $r$. Hence, using \eqref{eq2313-1} and \eqref{eq2313-2}, one deduces that
\[
	\mathscr P_\Omega(\{g>r\})\ge\sqrt{\mathscr{C}|\{g>r\}|}\ge\sqrt{\mathscr{C}s}
\]
for almost every $r\in(g^*(t),g^*(s))$. Integrating this inequality and using \eqref{eqn:1522-1} and \eqref{eqn:1451-1}, we obtain
\begin{align}\label{2142-1}
	\sqrt{\mathscr{C}s}\,[g^*(s)-g^*(t)]&\le\int_{g^*(t)}^{g^*(s)}\mathscr P_\Omega(\{g>r\})\,dr\nonumber\\
	&=\int_{\{g^*(t)<g<g^*(s)\}}|\nabla g|\,dx\nonumber\\
	&\le\|\nabla g\|_{L^\infty(\Omega)}|\{g^*(t)<g<g^*(s)\}|\nonumber\\
	&\le\|\nabla g\|_{L^\infty(\Omega)}(t-s).
\end{align}
Consequently, $g^*$ is Lipschitz on each $[\delta,S]$, $0<\delta<S$. Note that $g^*(S)=0$, since $|\{g>0\}|=S$. Taking $t=S$ in \eqref{2142-1}, one deduces
\[
	0\le g^*(s)\le\frac{\|\nabla g\|_{L^\infty(\Omega)}}{\sqrt{\mathscr{C}s}}(S-s),\qquad 0<s<S,
\]
and hence $\lim_{s\to S^-}g^*(s)=0$. Taking $t\searrow s$ in \eqref{2142-1}, we have
\begin{align}\label{eq2248-1}
	\sqrt{s}\,|(g^*)'(s)|\le\frac{1}{\sqrt{\mathscr{C}}}\|\nabla g\|_{L^\infty(\Omega)}\quad\text{for a.e. }s\in(0,S),
\end{align}
so that
\[
\int_0^S s|(g^*)'(s)|^2\,ds\le\frac{S}{\mathscr{C}}\|\nabla g\|_{L^\infty(\Omega)}^2<\infty.
\]

Let $R:=\sqrt{\frac{S}{\pi}}$. For $0<r<R$, set $F(r):=g^*(\pi r^2)$. According to \cite[Corollary~3.4.2, pp.~85--86]{Rak25}, we know that $g^*$ is locally absolutely continuous on $(0,|\Omega|)$. Thus, for $0<r_1<r_2<R$, we infer from \eqref{eq2248-1} that
\begin{align*}
	|F(r_2)-F(r_1)|&\le \int_{\pi r_1^2}^{\pi r_2^2}|(g^*)'(s)|\,ds\\
	&\le \frac{1}{\sqrt{\mathscr{C}}}\|\nabla g\|_{L^\infty(\Omega)}\int_{\pi r_1^2}^{\pi r_2^2}s^{-\frac{1}{2}}\,ds\\
	&=\frac{2\sqrt{\pi}}{\sqrt{\mathscr{C}}}\|\nabla g\|_{L^\infty(\Omega)}(r_2-r_1).
\end{align*}
Hence $F$ is Lipschitz on $(0,R)$ and has finite limits at both endpoints. Since $\lim_{s\to S^-}g^*(s)=0$, we have $\lim_{r\to R^-}F(r)=0$. Define $F(0):=\lim_{r\downarrow0}F(r)$ and extend $F$ by zero to $[R,\infty)$. The resulting function is Lipschitz on $[0,\infty)$ with the same Lipschitz bound.

For all $x,y\in\mathbb R^2$, we consequently have
\begin{align*}
	|F(|x|)-F(|y|)|&\le \frac{2\sqrt{\pi}}{\sqrt{\mathscr{C}}}\|\nabla g\|_{L^\infty(\Omega)}\bigl||x|-|y|\bigr|\\
	&\le \frac{2\sqrt{\pi}}{\sqrt{\mathscr{C}}}\|\nabla g\|_{L^\infty(\Omega)}|x-y|.
\end{align*}
Thus $x\mapsto F(|x|)$ is globally Lipschitz, is supported in $\overline{B_R}$, and agrees with $G$ almost everywhere. Furthermore, due to \cite[Theorem 4.5]{EG15}, we see that $G$ belongs to $W^{1,2}(\mathbb{R}^2)$. Finally, the chain rule and polar coordinates imply \eqref{a.4}.
\end{proof}

\begingroup
\interlinepenalty=10000
\footnotesize\linespread{1}\selectfont
\setlength{\bibsep}{0.5pt}
\renewcommand{\bibsection}{%
	\section*{\refname}%
	\addcontentsline{toc}{section}{\refname}%
}

\endgroup
\end{document}